\documentclass[12pt]{article}
\usepackage{latexsym,amsmath,amscd,amssymb,graphics}
\usepackage{enumerate}

\newcommand{\bfi}{\bfseries\itshape}

\makeatletter

\@addtoreset{figure}{section}
\def\thefigure{\thesection.\@arabic\c@figure}
\def\fps@figure{h, t}
\@addtoreset{table}{bsection}
\def\thetable{\thesection.\@arabic\c@table}
\def\fps@table{h, t}
\@addtoreset{equation}{section}

\makeatother

\newcommand{\todo}[1]
  {\vspace{5 mm}\par \noindent \marginpar{\textsc{ToDo}} \framebox{\begin
  {minipage}[c]{0.95 \textwidth} \tt #1
\end{minipage}}\vspace{5 mm}\par}

\begin{document}

\newtheorem{theorem}{Theorem}[section]
\newtheorem{definition}[theorem]{Definition}
\newtheorem{lemma}[theorem]{Lemma}
\newtheorem{remark}[theorem]{Remark}
\newtheorem{proposition}[theorem]{Proposition}
\newtheorem{corollary}[theorem]{Corollary}
\newtheorem{example}[theorem]{Example}

\makeatletter

\title{Extensions of Banach Lie-Poisson spaces}
\author{Anatol Odzijewicz$^{1}$ and Tudor S. Ratiu$^{2}$}
\addtocounter{footnote}{1}
\footnotetext{Institute of Physics, University of Bialystok, Lipowa 41, PL-15424
Bialystok, Poland. \texttt{aodzijew@labfiz.uwb.edu.pl}}
\addtocounter{footnote}{1}
\footnotetext{Centre Bernoulli,
\'Ecole Polytechnique F\'ed\'erale de Lausanne. CH--1015 Lausanne.
Switzerland. \texttt{Tudor.Ratiu@epfl.ch} }
\date{September 9, 2003}
\maketitle

\makeatother

\maketitle

%|||-------------------text width----------------------|||

\noindent \textbf{AMS Classification:}
 46L10, 46T05, 46T20, 53D17, 53Z05

\begin{abstract}
The extension of Banach Lie-Poisson spaces is studied and linked to
the extension of a special class of Banach Lie algebras. The
 case of $W ^\ast$-algebras is given particular attention.
Semidirect products and the extension of the restricted
Banach Lie-Poisson space by the Banach Lie-Poisson space of compact
operators are given as examples.
\end{abstract}

\tableofcontents
%-------------------------------------------------------
%--------------------------------------------------------
\section{Introduction}
\label{Section: Introduction}

The dual of any finite dimensional Lie algebra carries
a linear Poisson bracket, called \textit{Lie-Poisson structure\/},
which is pervasive in classical mechanics. Many Hamiltonian
systems, such as the free or heavy rigid body equations, the finite
and periodic Toda lattice, the geodesics on quadrics, or the
Neumann and Rosochatius system, have alternate non-canonical
descriptions in Lie-Poisson formulation. Formally, several
evolutionary partial differential equations, such as the ideal non
viscous fluid, ideal magnetohydrodynamics, the Poisson-Vlasov,
Korteweg-de Vries, Kadomtsev-Petviashvili, or the linear and
nonlinear wave and  Schr\"odinger equations also have Lie-Poisson
formulations (see, for example, Marsden and Ratiu [1994], Dubrovin,
Krichever, and Novikov [2001] and references therein). While these
applications of this linear Poisson structure emerged in the work of
the past decades, the structure itself goes back to Lie [1890], who
introduced it simultaneously with the concept of Lie algebra.
The theory of Lie-Poisson spaces in finite dimensions is
complete and is part of Poisson
geometry (see, for example, Weinstein [1983], [1998], Marsden and
Ratiu [1994], Vaisman [1994]). In infinite dimensions a general
theory is lacking. Chernoff and Marsden [1974] represents the first
systematic attempt to lay the foundations of Banach symplectic
geometry.

Motivated by our understanding of infinite dimensional Hamiltonian
systems, questions surrounding the notion of momentum map, and
problems in quantum mechanics including the theory of coherent
states,  Odzijewicz and Ratiu [2003] propose a definition of Banach
Lie-Poisson space and link it to classical and quantum reduction,
the theory of $W^\ast$-algebras, and momentum maps in infinite
dimensions. Banach Lie-Poisson spaces naturally appear in this
context as preduals of Banach Lie algebras. For example, the spaces
of compact and trace class operators on a complex separable Hilbert
space carry a natural Lie-Poisson bracket. The present work
develops further this point of view by addressing the fundamental
question of construction of new Banach Lie-Poisson spaces out of
given ones. One such scheme is given by the method of extensions.

The problem of extension in various categories plays a central role
in the understanding  of its objects and morphisms. It gives a
method to construct new objects out of old ones whose properties
are then well understood. The category of Banach Lie-Poisson
spaces is no
exception. The goal of this paper is to present the theory of
extensions for Banach Lie-Poisson spaces and to give several
mathematically and physically relevant examples.

The paper is organized as follows. In \S  \ref{Section: Banach
Lie-Poisson spaces} the necessary information on Banach Lie-Poisson
spaces found in Odzijewicz and Ratiu [2003] is collected. Only some
definitions and theorems necessary for the
subsequent development are given. With this background, the theory
of exact sequences of Banach Lie-Poisson spaces is presented in \S
\ref{Section: Exact sequences of Banach spaces}. It is shown that
exactness in this category is equivalent to exactness of the dual
Banach Lie algebra sequence in the subcategory of Banach Lie
algebras admitting a predual. Special attention is devoted to the
important case of the Banach Lie-Poisson spaces that are preduals
of $W ^\ast$-algebras. It is shown that if the the dual sequence is
exact in the category of
$W ^\ast$-algebras, then the exact sequence is necessarily that of
a direct sum of Banach Lie-Poisson spaces. Extensions of Banach Lie
algebras are discussed in \S \ref{Section: Extensions of Banach Lie
algebras}. All possible brackets on a Banach space direct sum of
Banach Lie algebras are characterized. With this preparation, \S
\ref{Section: Extensions of Banach Lie-Poisson spaces} presents all
extensions of Banach Lie-Poisson spaces underlying a Banach space
direct sum. Semidirect products of Banach Lie-Poisson
spaces with cocycles are a particular case of this theory. Even
more special, the case of the predual of the semidirect product of
a $W ^\ast$-algebra with a representation is treated in detail. The
example of the extension of the restricted Banach Lie-Poisson
space by the space of compact operators, important in
the theory of loop groups (Pressley and
Segal [1986], Wurzbacher [2001]), is also worked out.

\section{Banach Lie algebras and Lie-Poisson spaces}
\label{Section: Banach Lie-Poisson spaces}

This section briefly reviews the minimal background from
Odzijewicz and Ratiu [2003] necessary for the rest of
the paper. No proofs will be given here since they can
be found in the aforementioned paper.

Given a Banach
space $\mathfrak{b}$, the notation $\mathfrak{b}^\ast$ will always
be used for the Banach space dual to $\mathfrak{b}$. For $x \in
\mathfrak{b}^\ast$ and $b \in \mathfrak{b}$, the notation $\langle
x, b \rangle$ means the value of $x$ on $b$. Thus $\langle
\cdot, \cdot \rangle : \mathfrak{b}^\ast \times \mathfrak{b}
\rightarrow \mathbb{R}$ will  denote the natural bilinear continuous
duality pairing between $\mathfrak{b}$ and its dual
$\mathfrak{b}^\ast$. The notation $\mathfrak{b}_\ast$ will be
reserved for the {\bfi predual\/} of $\mathfrak{b}$, that
is, $\mathfrak{b}_\ast$ is a Banach space whose dual is
$\mathfrak{b}$. The predual is not unique, in general.
Note also that $\mathfrak{b}_\ast \hookrightarrow
\mathfrak{b}^\ast$ canonically and that $\mathfrak{b}_\ast $ is a
closed subspace of $\mathfrak{b}^\ast $.

Recall that a {\bfi Banach Lie algebra} $(\mathfrak{g},
[\cdot, \cdot])$ is a Banach space that is also a Lie
algebra such that the Lie bracket is a bilinear
continuous map $\mathfrak{g} \times \mathfrak{g}
\rightarrow \mathfrak{g}$. Thus the adjoint and
coadjoint maps $\operatorname{ad}_x :\mathfrak{g}
\rightarrow \mathfrak{g}$, $\operatorname{ad}_x y:=
[x,y]$, and $\operatorname{ad}_x^\ast : \mathfrak{g}^\ast
\rightarrow \mathfrak{g}^\ast$ are also continuous for
each $x \in \mathfrak{g}$.

A {\bfi Banach Poisson manifold\/} is a pair $(P, \{\cdot, \cdot
\})$ consisting of a smooth (real or complex) Banach manifold $P $
and a bilinear operation $\{\cdot, \cdot \}$ on the ring
$C^\infty (P)$ such that:
\begin{itemize}
\item $(C^\infty (P), \{\cdot, \cdot \})$ is Lie algebra,
\item the Leibniz identity holds: $\{fg, h\} = f\{g, h\} + \{f, h
\} g $ for all $f, g, h \in C^\infty (P) $,
\item for each $f \in C^\infty (P)$, the derivation $X_f:=
\{\cdot , f\}$ which is, in general, a section of $T^{\ast \ast}P$,
is a vector field on $P $.
\end{itemize}

A {\bfi Banach Lie-Poisson space\/}
$(\mathfrak{b},\{\cdot, \cdot\})$ is defined to be a
(real or complex) Banach space $\mathfrak{b}$ that is also a Poisson
manifold satisfying the additional condition that its dual
$\mathfrak{b}^\ast \subset C^\infty(\mathfrak{b})$ is a Banach Lie
algebra under the Poisson bracket operation. The following
characterization is crucial throughout this paper.
\begin{theorem}
\label{general theorem}
The Banach space $\mathfrak{b}$ is a Banach Lie-Poisson
space $(\mathfrak{b}, \{\cdot, \cdot \})$ if and only if
its dual $\mathfrak{b}^\ast$ is a Banach Lie algebra
$(\mathfrak{b}^\ast, [\cdot, \cdot ])$ satisfying
$\operatorname{ad}_x^\ast \mathfrak{b} \subset
\mathfrak{b} \subset \mathfrak{b}^{\ast \ast}$ for all
$x \in \mathfrak{b}^\ast$. Moreover, the Poisson bracket
of $f, g \in C^\infty(\mathfrak{b})$ is given by
\begin{equation}
\label{general LP}
\{f, g\}(b) = \langle [Df(b), Dg(b)], b \rangle,
\end{equation}
where $b \in \mathfrak{b}$ and $D$ denotes the Fr\'echet
derivative. If $h \in C^\infty(\mathfrak{b})$, the
associated Hamiltonian vector field is given by
\begin{equation}
\label{general Hamiltonian vector field}
X_h(b) = - \operatorname{ad}^\ast _{Dh(b)} b.
\end{equation}
\end{theorem}

A {\bfi morphism\/} between two
Banach Lie-Poisson spaces $\mathfrak{b}_1$ and $\mathfrak{b}_2$ is
a continuous linear map $\phi:\mathfrak{b}_1 \rightarrow
\mathfrak{b}_2$ that preserves the Poisson bracket, that is,
\[
\{f \circ \phi, g\circ \phi \}_1 = \{f, g \}_2 \circ \phi
\]
for any $f, g \in C^\infty(\mathfrak{b}_2)$. Such a map $\phi$
is  also called a {\bfi linear Poisson map}.
Define the {\bfi category\/} $\mathfrak{B}$ {\bfi of Banach
Lie-Poisson spaces\/} as the category whose objects are the Banach
Lie-Poisson spaces and whose morphisms are the linear Poisson maps.

Let $\mathfrak{L}$ denote the {\bfi category of Banach Lie
algebras}; its objects are Banach Lie algebras and its morphisms
are continuous Lie algebra homomorphisms.

Denote by $\mathfrak{L}_0$ the following subcategory of
$\mathfrak{L}$. An object of
$\mathfrak{L}_0$ is a Banach Lie algebra $\mathfrak{g}$ admitting
a predual $\mathfrak{g}_\ast$, that is, $(\mathfrak{g}_\ast)^\ast
= \mathfrak{g}$, and satisfying $\operatorname{ad}^\ast_{\mathfrak{g}}
\mathfrak{g}_\ast \subset \mathfrak{g}_\ast$ where $\operatorname{ad}^\ast$ is the
coadjoint  representation of $\mathfrak{g}$ on $\mathfrak{g}^\ast$;
recall that
$\mathfrak{g}_\ast$ is a closed subspace of $\mathfrak{g}^\ast$. A
morphism in the category $\mathfrak{L}_0$ is a Banach Lie algebra
homomorphism $\psi: \mathfrak{g}_1 \rightarrow \mathfrak{g}_2$ such
that the dual map $\psi^\ast: \mathfrak{g}_2^\ast \rightarrow
\mathfrak{g}_1^\ast$ preserves at least one choice of the
corresponding preduals, that is, $\psi^\ast: (\mathfrak{g}_2)_\ast
\rightarrow (\mathfrak{g}_1)_\ast$, where $(\mathfrak{g}_i)_\ast$
is one possible predual of $\mathfrak{g}_i$ for $i = 1,2 $. Let
$\mathfrak{L}_{0u}$ be the subcategory  of
$\mathfrak{L}_0$ whose objects have a
\textit{unique\/} predual.

\begin{theorem}
\label{functor}
There is a contravariant functor $\mathfrak{F}: \mathfrak{B} \rightarrow
\mathfrak{L}_0$ defined by $\mathfrak{F}(\mathfrak{b}) = \mathfrak{b}^\ast$ and
$\mathfrak{F}(\phi) = \phi^\ast$. On the subcategory
$\mathfrak{F}^{-1}(\mathfrak{L}_{0u}) \subset \mathfrak{B}$ this
functor is  invertible. The inverse of $\mathfrak{F}$ is given by
$\mathfrak{F}^{-1}(\mathfrak{g}) = \mathfrak{g}_\ast$ and
$\mathfrak{F}^{-1}(\psi) = \psi^\ast|_{(\mathfrak{g}_2)_\ast}$,
where $\psi: \mathfrak{g}_1 \rightarrow \mathfrak{g}_2$.
\end{theorem}

\medskip

The internal structure of morphisms in $\mathfrak{B}$ is given by
the following results.

\begin{proposition}
\label{decomposition}
Let $\phi:\mathfrak{b}_1 \rightarrow
\mathfrak{b}_2$ be a linear Poisson map between Banach Lie-Poisson
spaces and assume that $\operatorname{im} \phi$ is closed in
$\mathfrak{b}_2$. Then the Banach space $\mathfrak{b}_1/ \ker
\phi$ is predual to $\mathfrak{b}_2^\ast/ \ker \phi^\ast$, that
is, $(\mathfrak{b}_1/ \ker \phi)^\ast \cong \mathfrak{b}_2^\ast/
\ker \phi^\ast$. In addition, $\mathfrak{b}_2^\ast/ \ker
\phi^\ast$ is a Banach Lie algebra satisfying the condition
$\operatorname{ad}_{[x]}^\ast \left(\mathfrak{b}_1/ \ker \phi
\right) \subset \mathfrak{b}_1/ \ker \phi$ for all $[x] \in
\mathfrak{b}_2^\ast/ \ker \phi^\ast$ and $\mathfrak{b}_1/\ker
\phi$ is a Banach Lie-Poisson space. Moreover, the following
properties hold:
\begin{enumerate}
\item[{\rm (i)}] the quotient map $\pi: \mathfrak{b}_1
\rightarrow \mathfrak{b}_1/\ker \phi$ is a surjective
linear Poisson map;
\item[{\rm (ii)}] the map $\iota: \mathfrak{b}_1/\ker \phi
\rightarrow \mathfrak{b}_2$ defined by
$\iota([b]):= \phi(b)$, where $b \in \mathfrak{b}_1$ and
$[b] \in \mathfrak{b}_1/\ker \phi$ is an injective linear Poisson map;
\item[{\rm (iii)}] the decomposition $\phi= \iota \circ \pi$
into a surjective and an injective linear Poisson
map is valid.
\end{enumerate}
\end{proposition}
Proposition \ref{decomposition} reduces the study of
linear Poisson maps with closed range between Banach Lie-Poisson
spaces to the study of surjective  and injective linear Poisson
maps, which is carried out in the next propositions.

\begin{proposition}
\label{surjective LP map}
Let $(\mathfrak{b}_1, \{\cdot, \cdot \})$ be a Banach Lie-Poisson space
and let $\pi : \mathfrak{b}_1 \rightarrow \mathfrak{b}_2$ be a continuous linear
surjective map onto the Banach space $\mathfrak{b}_2$. Then $\mathfrak{b}_2$
carries a Banach Lie-Poisson structure such that $\pi$ is a
linear Poisson map  if and only if
$\operatorname{im} \pi^\ast \subset \mathfrak{b}_1^\ast$ is closed
under the Lie bracket $[\cdot, \cdot]_1$ of $\mathfrak{b}_1^\ast$.
This Banach Lie-Poisson structure on $\mathfrak{b}_2 $ is unique
and it is called the {\bfi coinduced\/} structure by the mapping
$\pi$. The map $\pi^\ast :\mathfrak{b}_2^\ast \rightarrow
\mathfrak{b}_1^\ast$ is a Banach Lie algebra morphism whose dual
$\pi^{\ast \ast}:
\mathfrak{b}_1^{\ast \ast} \rightarrow \mathfrak{b}_2^{\ast \ast}$
maps $\mathfrak{b}_1$ into $\mathfrak{b}_2$.
\end{proposition}

\begin{proposition}
\label{injective LP map}
Let $\mathfrak{b}_1$ be a Banach space,
$(\mathfrak{b}_2, \{\cdot, \cdot \}_2)$ be a Banach Lie-Poisson
space, and $\iota : \mathfrak{b}_1 \rightarrow \mathfrak{b}_2$ be
an injective continuous linear map with closed range. Then
$\mathfrak{b}_1$ carries a unique Banach Lie-Poisson structure such that
$\iota$ is a linear Poisson map if and only if
$\ker \iota^\ast$ is an ideal in the Banach Lie algebra
$\mathfrak{b}_2^\ast$. This Banach Lie-Poisson structure on
$\mathfrak{b}_1 $ is unique and it is called the  structure {\bfi
induced\/} by the mapping $\iota$. The
map $\iota^\ast :\mathfrak{b}_2^\ast \rightarrow
\mathfrak{b}_1^\ast$ is a Banach Lie algebra morphism whose dual
$\iota^{\ast \ast}:
\mathfrak{b}_1^{\ast \ast} \rightarrow \mathfrak{b}_2^{\ast \ast}$
maps $\mathfrak{b}_1$ into $\mathfrak{b}_2$.
\end{proposition}

For later applications we shall need also the notion of the product
of  Banach Poisson manifolds.

\begin{theorem}
\label{product theorem}
Given the Banach Poisson manifolds $(P_1, \{\,,\}_1)$ and
$(P_2, \{\,,\}_2)$ there is a unique Banach Poisson structure
$\{\,,\}_{12}$ on the product manifold $P_1 \times P_2$ such that:
\begin{enumerate}
\item[{\rm (i)}] the canonical projections
$\pi_1: P_1 \times P_2 \rightarrow P_1$
and $\pi_2: P_1 \times P_2 \rightarrow P_2$ are Poisson maps;
\item[{\rm (ii)}] $\pi_1^\ast (C^\infty(P_1))$ and $\pi_2^\ast (C^\infty(P_2))$
are Poisson commuting subalgebras of $C^\infty(P_1 \times P_2)$.
\end{enumerate}
This unique Poisson structure on $P_1 \times P_2$ is called the
{\bfi product\/} Poisson structure and its bracket is given by the formula
\begin{equation}
\label{product bracket}
\{f, g\}_{12}(p_1, p_2) = \{f_{p_2}, g_{p_2}\}_1(p_1) +
\{f_{p_1}, g_{p_1}\}_2(p_2),
\end{equation}
where $f_{p_1}, g_{p_1}  \in C^\infty(P_2)$ and
$f_{p_2}, g_{p_2} \in C^\infty(P_1)$ are the partial functions given by
$f_{p_1}(p_2) = f_{p_2}(p_1) = f(p_1, p_2)$ and similarly for $g$.
\end{theorem}

\bigskip

\section{Exact sequences of Banach Lie-Poisson spaces}
\label{Section: Exact sequences of Banach spaces}

In this section we will study exact sequences in the
categories presented in  the previous section.
Exactness in the categories $ \mathfrak{L} $, $\mathfrak{L}_0 $,
$\mathfrak{L}_{0u} $, and
$\mathfrak{B}$ is defined in the following way.

\begin{definition}
\label{def: exact sequence of lie algebras}
A sequence of Banach Lie algebras
\begin{equation}
\label{dual lie algebra sequence}
0 \longrightarrow \mathfrak{n} \stackrel{\iota}\longrightarrow
\mathfrak{g} \stackrel{\pi}\longrightarrow \mathfrak{h}
\longrightarrow 0
\end{equation}
is {\bfi exact\/} if it is exact as a sequence in the category of
Banach spaces and all maps are Banach Lie algebra homomorphisms.
The Lie algebra $\mathfrak{g}$ is said to be an {\bfi
extension} of $\mathfrak{h}$ by $\mathfrak{n}$.
\end{definition}

\begin{definition}
\label{def: exact sequence in l zero}
In the categories $\mathfrak{L}_0 $ (respectively
$\mathfrak{L}_{0u}$) the  sequence \eqref{def: exact sequence of lie
algebras} is {\bfi exact\/} if it  is exact in the category
$\mathfrak{L}$ and the duals of the maps in the sequence preserve
at least one choice of (respectively the uniquely associated)
predual spaces, that is $\iota(\mathfrak{g}_\ast ) \subset
\mathfrak{n}_\ast $ and $\pi_\ast (\mathfrak{h}_\ast) \subset
\mathfrak{g}_\ast $, where $\mathfrak{n}_\ast, \mathfrak{g}_\ast,
\mathfrak{h}_\ast $ are preduals of $\mathfrak{n},
\mathfrak{g}, \mathfrak{h}$ respectively. Like in the previous case,
$\mathfrak{g}$ is said to be an {\bfi extension} of $\mathfrak{h}$
by $\mathfrak{n}$ in these two categories.
\end{definition}

\begin{definition}
\label{def: exact sequence of lie poisson spaces}
A sequence of Banach Lie-Poisson spaces
\begin{equation}
\label{lie poisson sequence}
0 \longrightarrow \mathfrak{a} \stackrel{j}\longrightarrow
\mathfrak{b} \stackrel{p}\longrightarrow \mathfrak{c}
\longrightarrow 0
\end{equation}
is {\bfi exact\/} if it is exact as
a sequence in the category of Banach spaces and all maps are
linear Poisson maps.  The Banach Lie-Poisson space $\mathfrak{b}$ is
said to be an {\bfi extension\/} of $\mathfrak{c}$ by
$\mathfrak{a}$.
\end{definition}

The goal of this section is to study under what  conditions the
functor $\mathfrak{F}$ preserves exactness.  The answer is given by
the following theorem.

\begin{theorem}
\label{thm: exact sequence of lie poisson}
The Banach spaces $\mathfrak{a}$, $\mathfrak{b}$, $\mathfrak{c}$
form an exact sequence \eqref{lie poisson sequence}
of Banach Lie-Poisson spaces if and only if their duals
$\mathfrak{n} : = \mathfrak{c}^\ast$, $\mathfrak{g}: =
\mathfrak{b}^\ast$, $\mathfrak{h}: = \mathfrak{a}^\ast$ form an
exact sequence of Banach Lie algebras \eqref{dual lie algebra
sequence}  in the category $\mathfrak{L}_0$, where $\iota: = p
^\ast$ and $\pi: =  j ^\ast$.
In particular, if $\mathfrak{g}$ is the direct sum $\mathfrak{g}
\cong \mathfrak{n} \oplus \mathfrak{h}$ of Banach Lie algebras with
$\iota$ the inclusion of the first component and $\pi$
the projection on the second component, then
$\mathfrak{b}$ can be chosen as the direct sum
$\mathfrak{a}
\oplus\mathfrak{c}$ of the Banach Lie-Poisson spaces
$\mathfrak{a}$ and $\mathfrak{c}$ with $j $ the inclusion
on the first component and $p $ the  projection  on the second
component.
\end{theorem}

For the proof we shall need a few preparatory lemmas. The subjects of two first lemmas are well known facts but we include
them for the sake of completeness.

\begin{lemma}
\label{lemma: direct sum}
Let $U$ and $W$ be Banach spaces.  Then one
has the canonical isomorphism
\[
\left(U\oplus W\right)^\ast \cong
U^\ast \oplus W^\ast.
\]
\end{lemma}

\noindent \textbf{Proof.} To $f \in \left(U\oplus
W\right)^\ast $ associate the pair $(f|_U,
f|_W) \in U^\ast \oplus W^\ast$.
This map is clearly linear and continuous if on the direct sum one
takes the norm given by the sum of norms in each component. The map
that associates to $(\gamma, \alpha) \in U^\ast \oplus
W^\ast$ the functional $\gamma + \alpha \in \left(
U\oplus W\right)^\ast$, defined by $(\gamma + \alpha)(c, a) : =
\gamma(c) + \alpha(a) $, is also linear and continuous. The two
maps are clearly inverses of each other.
\quad $\blacksquare $

\begin{lemma}
\label{prop: exact dual sequence}
If one has the exact sequence of Banach spaces
\begin{equation}
\label{exact sequence}
0 \longrightarrow U \stackrel{\iota}\longrightarrow
V \stackrel{\pi}\longrightarrow W
\longrightarrow 0
\end{equation}
then the dual sequence
\begin{equation}
\label{dual exact sequence}
0 \longrightarrow W^\ast \stackrel{\pi^\ast}\longrightarrow
V^ \ast \stackrel{\iota^\ast}\longrightarrow U^\ast
\longrightarrow 0
\end{equation}
is also exact.
\end{lemma}

\noindent\textbf{Proof}. The linear continuous map $\pi ^\ast$ is
injective. Indeed, if $\gamma \in \ker \pi ^\ast \subset W ^\ast$,
then $( \gamma\circ \pi)(v) = 0 $ for all $v \in V $. Surjectivity
of $\pi$ implies then that $\gamma(w) = 0 $ for all $w \in W $ so
$\gamma = 0 $.

The linear continuous map $\iota ^\ast$ is surjective. Indeed,
since $\iota(U) = \ker \pi$ is a closed subspace of $V $ and
$\iota$ is injective,  $\iota: U
\rightarrow \iota(U) $ is a Banach space isomorphism. Then, if
$\alpha \in U ^\ast$, the linear functional $\alpha\circ \iota ^{-1}
: \iota(U) \rightarrow \mathbb{C}$ is continuous. Extend
this functional to $\beta\in V ^\ast$ by the Hahn-Banach Theorem.
Thus, for any $u \in U $, we have $\iota^\ast(\beta)(u) =
\beta(\iota(u)) = ( \alpha \circ \iota^{-1})(\iota(u)) =
\alpha(u)$, which shows that $\iota^\ast(\beta) = \alpha$, that
is, $\iota ^\ast$ is onto.

Since $\iota ^\ast \circ \pi ^\ast = ( \pi \circ \iota) ^\ast = 0
$ by exactness of the sequence \eqref{exact sequence}, it follows
that $\pi^\ast(W ^\ast ) \subset \ker \iota^\ast$. To prove the
opposite inclusion, let $\beta \in V ^\ast$ be such that
$\iota^\ast(\beta) =0 $. Define $\tilde \beta : W
\rightarrow \mathbb{C}$ by $\tilde \beta(w) =
\beta(v)$, if $w = \pi(v) $; thus $\tilde \beta\circ \pi=
\beta$. Since
$\beta|_{\iota(U)} = 0
$ by hypothesis and $\iota(U) = \ker \pi$ by exactness of
\eqref{exact sequence}, it follows that $\tilde \beta$ is well
defined. It is straightforward to verify that $\tilde \beta $ is
linear and continuous using the Banach space isomorphism
$V/\iota(U) \cong W $. Finally, $\beta =
\tilde \beta \circ \pi= \pi ^\ast(\tilde \beta) \in \pi^\ast(W
^\ast)$. \quad $\blacksquare $

\begin{lemma}
\label{prop: exact predual sequence}
Assume that all Banach spaces in the exact sequence \eqref{exact
sequence} admit preduals, that is, there are Banach spaces
$U_{\ast} $, $V_{\ast} $, and $W_{\ast}$ such that $U =
(U_{\ast})^\ast $, $V = (V_{\ast})^\ast $, and $W = (W_{\ast})^\ast
$, respectively. Assume, in addition, that  the dual maps $\pi^\ast$
and $\iota^\ast$ preserve the predual spaces, that is,
$\pi^\ast(W_\ast) \subset V_\ast $ and $\iota^\ast (V_\ast) \subset
U_\ast $. Then one has the following commutative diagram of exact
sequences
\unitlength=5mm
\begin{center}
\begin{picture}(19,2)
%first line
\put(-3,0){\makebox(0,0){$0$}}
\put(1,0){\makebox(0,0){$W^\ast$}}
\put(9,0){\makebox(0,0){$V^\ast$}}
\put(16,0){\makebox(0,0){$U^\ast$}}
\put(21,0){\makebox(0,0){$0$}}
%second line
\put(-3,-5){\makebox(0,0){$0$}}
\put(1,-5){\makebox(0,0){$W_\ast$}}
\put(9,-5){\makebox(0,0){$V_\ast$}}
\put(16,-5){\makebox(0,0){$U_\ast$}}
\put(21,-5){\makebox(0,0){$0$}}
%third line
\put(1,-9){\makebox(0,0){$0$}}
\put(9,-9){\makebox(0,0){$0$}}
\put(16,-9){\makebox(0,0){$0$}}
%first vertical arrows
\put(1,-4){\vector(0,1){3}}
\put(1,-8){\vector(0,1){2}}
%second vertical arrows
\put(9,-4){\vector(0,1){3}}
\put(9,-8){\vector(0,1){2}}
%third vertical arrows
\put(16,-4){\vector(0,1){3}}
\put(16,-8){\vector(0,1){2}}
%first horizontal arrows
\put(-2.3,0){\vector(1,0){2.2}}
\put(2,0){\vector(1,0){6.4}}
\put(9.6,0){\vector(1,0){5.4}}
\put(16.8,0){\vector(1,0){3.5}}
%second horizontal arrows
\put(-2.4,-5){\vector(1,0){2.2}}
\put(2,-5){\vector(1,0){6.4}}
\put(9.6,-5){\vector(1,0){5.4}}
\put(16.8,-5){\vector(1,0){3.5}}
%first horizontal labels
\put(5,.7){\makebox(0,0){$\pi^\ast$}}
\put(12,.7){\makebox(0,0){$\iota^\ast$}}
%second horizontal labels
\put(5,-4.3){\makebox(0,0){$\pi_\ast$}}
\put(12,-4.3){\makebox(0,0){$\iota_\ast$}}
\end{picture}
\end{center}
\vspace{2in}
where all vertical arrows are inclusions and the maps in the
second line are defined by restriction, that is, $\pi_\ast : = \pi
^\ast|_{W_\ast}$ and $\iota_\ast : = \iota^\ast|_{V_\ast}$.

In particular, if $V \cong U \oplus W $, $U = (U_\ast)^\ast$,
$W = (W_\ast) ^\ast$, and the maps $\iota$ and $\pi$ in the
sequence
\[
0 \longrightarrow U \stackrel{\iota}\longrightarrow
U \oplus W \stackrel{\pi}\longrightarrow
W\longrightarrow 0
\]
are defined by $\iota(u) : = (u, 0)$, $\pi(u, w): = w $, for $u
\in U $, $w \in W $, then the commutative diagram above is valid
for $V_\ast = U_\ast \oplus W_\ast $.
\end{lemma}

\noindent\textbf{Proof}. Since $\pi^\ast(W_\ast) \subset V_\ast $
and $\iota^\ast (V_\ast) \subset U_\ast $, the maps in the second
line are well defined. The map $\pi_\ast $ is injective because it
is the restriction of the injective map $\pi ^\ast$ to the Banach
subspace $W_\ast \subset W^\ast $.

The relation $\pi\circ \iota = 0 $ implies that $\iota_\ast \circ
\pi_\ast = ( \iota ^\ast \circ \pi^\ast)|_{W_\ast} = ( \pi\circ
\iota) ^\ast |_{W_\ast} = 0 $ which shows that
\begin{equation}
\label{first inclusion at v lower star}
\operatorname{im} \pi_\ast \subset \ker \iota_\ast.
\end{equation}
Let us assume that there is some element $b \in \ker \iota_\ast $
such that $b \notin \overline{\operatorname{im}\pi_\ast} \subset
V_\ast$. Then, by the Hahn-Banach Theorem, there exists an element
$v \in V $ such that $v(b) \neq 0 $ and
$v|_{\overline{\operatorname{im}\pi_\ast}} \equiv 0 $. The
exactness of the sequence \eqref{exact sequence} and $b \in \ker
\iota_\ast \subset V_\ast \subset V^\ast $ implies that
\begin{equation}
\label{b on kernel of pi}
b|_{\ker \pi} = 0.
\end{equation}
However, one concludes from
$v|_{\overline{\operatorname{im}\pi_\ast}} \equiv 0 $ that
\[
0 = v(\pi_\ast (\gamma)) = \pi ^\ast(\gamma)(v) = \gamma(\pi(v))
\]
for any $\gamma\in W_\ast $, which implies that $\pi(v) = 0 $, that
is, $v \in \ker \pi$. By \eqref{b on kernel of pi} we have
$v(b) = 0 $, which contradicts the choice of $v $. Thus we have
proved that
\[
\overline{\operatorname{im}\pi_\ast} = \ker \iota_\ast.
\]
Thus, for any $b \in \ker \iota_\ast $ there is a sequence
$\{a_n\}_{n=1}^\infty
\subset W_\ast $ such that $\pi^\ast (a_n) = \pi_\ast (a_n)
\rightarrow b$ as $n \rightarrow \infty $ and, since
$\operatorname{im}\pi^\ast = \ker \iota^\ast $, there is some
element $a \in W^\ast $ such that $\pi^\ast(a) = b $. One obtains
form the above $\pi^\ast(a_n - a ) \rightarrow 0 $ as $n \rightarrow
\infty $. Now, because $\pi^\ast: W ^\ast \rightarrow \pi ^\ast(W
^\ast) = \ker \iota ^\ast $ is a continuous isomorphism
between Banach spaces, the Banach Isomorphism Theorem guarantees that
its inverse is also continuous and thus we have $a_n \rightarrow a $
as $n \rightarrow  \infty $. Since $W_\ast $ is closed in $W ^\ast$,
this implies that $a \in W_\ast $ and therefore $b \in
\operatorname{im} \pi _\ast $. This shows that
\begin{equation}
\label{exactness at v lower star}
\operatorname{im} \pi_\ast = \ker \iota_\ast.
\end{equation}

The last step of the proof is to show that
$\operatorname{im}\iota_\ast = U_\ast $. If
$\overline{\operatorname{im}\iota_\ast} \neq U_\ast $, then there
are elements $c \notin \overline{\operatorname{im}\iota_\ast} $
and $u \in U $ such that $u(c) \neq 0 $  and
$u|_{\overline{\operatorname{im}\iota_\ast}} \equiv 0 $. Thus one
has $0 = u(\iota_\ast (b)) = b(\iota(u)) $ for any $b \in V_\ast $,
which means that $\iota(u) = 0 $. Injectivity of $\iota$ implies
then that $u = 0 $, which contradicts the choice $u(c) \neq 0 $.
Thus we showed that
\[
\overline{\operatorname{im}\iota_\ast} = U_\ast.
\]
Therefore, for any $c \in U_\ast $ there is a sequence $\{ b_n
\}_{n=1}^\infty
\subset V_\ast $ and an element $b \in V^\ast $ such that
$\iota_\ast(b_n)
\rightarrow c = \iota^\ast(b) $ as $n \rightarrow \infty $. Thus
$\iota ^\ast(b_n - b) \rightarrow 0 $ as $n \rightarrow \infty $
and, using the isomorphism $V^\ast / \pi^\ast(W^\ast) \cong U^\ast
$, we conclude from the Banach Isomorphism Theorem that $[b_n - b]
\rightarrow 0 $ as $n \rightarrow \infty $ in   $V^\ast / \pi^\ast(W^\ast) $, where $[a]$ denotes the
equivalence class of $a \in V ^\ast$ in the quotient Banach
space $V^\ast / \pi^\ast(W^\ast)$. This means
that there  is a subsequence
$\{b_{n_k}\}_{k=1}^\infty
\subset V_\ast$ and an element $d \in W^\ast $ such that $b_{n_k} -
b
\rightarrow
\pi ^\ast(d) $ as $k \rightarrow \infty $. Therefore, since $V_\ast
$ is closed in $V^\ast $, we have that $\lim_{k \to \infty} b_{n_k}
\in V_\ast $ and one has
\[
c = \iota^\ast (b) = \iota ^\ast(\lim_{k \to \infty} b_{n_k} - \pi
^\ast(d)) = \iota^\ast(\lim_{k \to \infty} b_{n_k} )
\]
since $\iota ^\ast \circ \pi ^\ast = 0 $. This shows that $c \in
\operatorname{im} \iota_\ast $ which proves that $\operatorname{im}
\iota_\ast  = U_\ast $.

To prove the last statement, notice that by Lemma
\ref{lemma: direct sum} we have $(U_\ast \oplus W_\ast) ^\ast \cong
(U_\ast) ^\ast \oplus (W_\ast ) ^\ast \cong U \oplus W \cong V $,
that is,
$U_\ast \oplus W _\ast$ is a predual of $V$. In addition, it is
easy to see that if $f \in W_\ast \subset W ^\ast$ and $(g_1,
g_2 ) \in U_\ast \oplus W_\ast \subset U ^\ast\oplus W ^\ast$, then
$\pi^\ast(f) = (0, f) \in U _\ast \oplus W_\ast$ and
$\iota^\ast(g_1, g_2) = g_1 \in U_\ast$ and thus the hypotheses in
the first part of the lemma are verified.
\quad $\blacksquare $

\medskip

We have now the necessary background to prove the main theorem of
this section.
\medskip

\noindent\textbf{Proof of Theorem \ref{thm: exact sequence of lie
poisson}}.  Let us assume that $\mathfrak{a}, \mathfrak{b}$, and
$\mathfrak{c}$ form an exact sequence of Banach Lie-Poisson spaces
in the sense of Definition \ref{def: exact sequence of lie poisson
spaces}.  Lemma \ref{prop: exact dual sequence} guarantees that
their duals also form an exact sequence
\[
0 \longrightarrow \mathfrak{c}^\ast \stackrel{p^\ast}\longrightarrow
\mathfrak{b}^\ast \stackrel{j^\ast}\longrightarrow \mathfrak{a}^\ast
\longrightarrow 0
\]
of Banach spaces. From Propositions \ref{surjective LP map} and
\ref{injective LP map} it follows that this sequence is an exact
sequence of Banach Lie algebras in the category $\mathfrak{L}_0$.

Conversely, assume that $\mathfrak{n}, \mathfrak{g}$, and
$\mathfrak{h}$ form an exact sequence of Banach Lie algebras in the
category $\mathfrak{L}_0$ (see Definition \ref{def: exact sequence
in l zero}). Thus they have preduals $\mathfrak{n}_\ast,
\mathfrak{g}_\ast $, and $\mathfrak{h}_\ast $ which, by Lemma
\ref{prop: exact predual sequence}, form an exact sequence
\[
0 \longrightarrow \mathfrak{h}_\ast
\stackrel{\pi_\ast}\longrightarrow
\mathfrak{g}_\ast \stackrel{\iota_\ast}\longrightarrow
\mathfrak{n}_\ast
\longrightarrow 0
\]
of Banach spaces. Now, again by Propositions \ref{surjective LP
map},  \ref{injective LP map}, and Theorem \ref{general theorem},
the preduals $\mathfrak{n}_\ast,
\mathfrak{g}_\ast $, and $\mathfrak{h}_\ast $ are Banach Lie
Poisson spaces and the maps $\pi_\ast $ and $\iota_ \ast $ are
linear Poisson maps.

The last statement of the theorem is proved in the following way.
By the second part of Lemma \ref{prop: exact predual sequence},
$\mathfrak{a}\oplus\mathfrak{c}$ can be taken as the predual space
to $\mathfrak{h} \oplus \mathfrak{n}$. By the first part of the
theorem, the natural maps $j $ and $p $ are linear Poisson maps.
The desired conclusion now immediately follows from Theorem
\ref{product theorem}, if we show that the spaces $\pi_\mathfrak{a}
^\ast (C^\infty(\mathfrak{a}))$ and $\pi_\mathfrak{c}
^\ast (C^\infty(\mathfrak{c}))$
are Poisson commuting subalgebras of $C^\infty(\mathfrak{a} \oplus
\mathfrak{c})$, where $\pi_\mathfrak{a}: \mathfrak{a}
\oplus\mathfrak{c} \rightarrow \mathfrak{a}$ and $\pi_
\mathfrak{c}: \mathfrak{a} \oplus\mathfrak{c}\rightarrow
\mathfrak{c}$ are the natural projections. If $f \in
C^\infty(\mathfrak{a})$ and $g \in C^\infty(\mathfrak{c})$, then
$Df(a) \in \mathfrak{a}^\ast = \mathfrak{h}$ and $Dg(c) \in
\mathfrak{c}^\ast = \mathfrak{n}$, so that
$[D(\pi_\mathfrak{a}^\ast f)(a, c), D(\pi_\mathfrak{c}^\ast g)(a,
c)] = [(Df(a), 0), (0, Dg(c))] = 0$, since $\mathfrak{h} \oplus
\mathfrak{n}$ is a direct sum of Lie algebras. Formula
\eqref{general LP} of the Lie-Poisson bracket insures then that
$\{\pi_\mathfrak{a}^\ast f, \pi_\mathfrak{c}^\ast g \} = 0 $ as
required.
\quad $\blacksquare $
\bigskip

Therefore, one concludes from Theorem \ref{thm: exact sequence of lie
poisson} that the
problem of extensions in the category of Banach Lie-Poisson spaces
is equivalent to that in the category $\mathfrak{L}_0 $.

It was shown in Odzijewicz and Ratiu [2003] that on the predual
Banach space of a $W ^\ast$-algebra (von Neumann algebra) there is
a canonically defined Banach Lie-Poisson structure. Since the
theory of von Neumann algebras  is closely related to crucial
problems of quantum physics (Emch [1972], Bratteli [1979], [1981])
we shall discuss the application of Theorem \ref{thm: exact sequence
of lie poisson} to this subcase.

Recall that a $W^*$-algebra is a $C^*$-algebra
$\mathfrak{m}$ which posses a predual Banach space
$\mathfrak{m}_\ast$, i.e.
$\mathfrak{m}=(\mathfrak{m}_\ast)^*$; this predual is unique
(Sakai [1971]). Since
$\mathfrak{m}^\ast=(\mathfrak{m}_\ast)^{**}$, the predual Banach
space $\mathfrak{m}_\ast$ canonically embeds into the Banach space
$\mathfrak{m}^\ast$ dual to $\mathfrak{m}$. Thus we shall always
think of $\mathfrak{m}_\ast$ as a Banach subspace of
$\mathfrak{m}^\ast$. The existence of $\mathfrak{m}_\ast$ allows
the introduction of the
$\sigma(\mathfrak{m},\mathfrak{m}_\ast)$-topology on the
$W^*$-algebra
    $\mathfrak{m}$; for simplicity we shall call it the $\sigma$-topology in
the sequel. Recall that a net $\{x_\alpha\}_{\alpha \in A} \subset
\mathfrak{m}$ converges to $x \in \mathfrak{m}$ in the $\sigma$-topology
if, by definition, $\lim_{\alpha \in A} \langle x_\alpha, b \rangle =
\langle x, b \rangle$ for all $b \in \mathfrak{m}_\ast$.
The predual space $\mathfrak{m}_\ast$ is characterized as the
subspace of $\mathfrak{m}^\ast$ consisting of all
$\sigma$-continuous linear functionals on $\mathfrak{m}$ (Sakai
[1971]).

\medskip

A {\bfi homomorphism of $W ^\ast$-algebras\/} $\phi: \mathfrak{m}_1
\rightarrow \mathfrak{m}_2 $ is a $\sigma$-continuous $\ast
$-algebra homomorphism (and is hence automatically norm continuous;
Sakai [1971]). Note that $\phi ^\ast(\mathfrak{m}_{2\ast})
\subset \mathfrak{m}_{1\ast}$, where
$\mathfrak{m}_{i\ast}$ is the unique predual of
$\mathfrak{m}_i$, for $i = 1,2$. Indeed, since any element $b
\in \mathfrak{m}_{2\ast}$ is $\sigma$-continuous on
$\mathfrak{m}_2 $ and $\phi$ is also
$\sigma$-continuous, it follows that $\phi ^\ast(b) = b \circ
\phi$ is $\sigma$-continuous on $\mathfrak{m}_1 $ and hence is an
element of $\mathfrak{m}_{1\ast}$.

Conversely, assume that $\mathfrak{m}_1 $ and $\mathfrak{m}_2
$ are $W ^\ast$-algebras and that $\phi: \mathfrak{m}_1 \rightarrow
\mathfrak{m}_2 $ is  a $\ast $-homomorphism satisfying
$\phi ^\ast(\mathfrak{m}_{2\ast})\subset \mathfrak{m}_{1\ast}$.
Then $\phi$ is $\sigma$-continuous. Indeed, if
$\{x_\alpha\}_{\alpha\in A } \subset \mathfrak{m}_1$ is a net
$\sigma$-converging to $x \in \mathfrak{m}_1 $, then for any
$b_1\in \mathfrak{m}_{1\ast}$ we have $\lim_{\alpha \in A} \langle
x_\alpha, b_1 \rangle = \langle x, b_1 \rangle$. We have for any
$b_2 \in \mathfrak{m}_{2\ast}$,
$\lim_{\alpha \in A} \langle \phi(x_\alpha), b_2 \rangle =
\lim_{\alpha \in A} \langle x_\alpha, b_2\circ \phi \rangle =
\langle x, b_2\circ \phi \rangle$, since, by hypothesis, $b_2 \circ
\phi = \phi ^\ast(b_2) \in \mathfrak{m}_{1\ast}$. This shows that
$\lim_{\alpha \in A} \langle \phi(x_\alpha), b_2 \rangle = \langle
\phi(x), b_2 \rangle $ for any $b_2 \in \mathfrak{m}_{2 \ast}$,
that is, $\phi: \mathfrak{m}_1 \rightarrow \mathfrak{m}_2 $ is
$\sigma$-continuous. These arguments prove the following.
\begin{proposition}
\label{w star homomorphism condition}
Let $\mathfrak{m}_1 $ and $\mathfrak{m}_2 $ be $W ^\ast$-algebras
and $\phi: \mathfrak{m}_1 \rightarrow \mathfrak{m}_2 $ a $\ast
$-homomorphism. Then $\phi$ is  a $W ^\ast$-algebra homomorphism if
and only if $\phi^\ast$ preserves the preduals, that is, $\phi
^\ast(\mathfrak{m}_{2\ast})\subset
\mathfrak{m}_{1\ast}$.
\end{proposition}

Since any $W ^\ast$-algebra is a Banach Lie algebra relative to
the commutator bracket and possesses a unique predual, this
proposition plus the condition $\operatorname{ad}_a^\ast
\mathfrak{m}_\ast \subset \mathfrak{m}_\ast \subset
\mathfrak{m}^\ast $ for any $a \in \mathfrak{m}$, shows that
$\mathfrak{W}$ is a subcategory of $\mathfrak{L}_{0u}$ (see Theorem
\ref{general theorem}). The condition $\operatorname{ad}_a^\ast
\mathfrak{m}_\ast \subset \mathfrak{m}_\ast$ for all $a \in
\mathfrak{m}_\ast $ is always satisfied. Indeed,  left and right
multiplication by $a\in\mathfrak{m}$ define uniformly and
$\sigma$-continuous maps $L_a: \mathfrak{m}\ni x\mapsto
ax\in\mathfrak{m}$ and $R_a: \mathfrak{m} \ni x\mapsto
xa\in\mathfrak{m}$ (Sakai [1971]). Let
$L_a^*:\mathfrak{m}^\ast\to\mathfrak{m}^\ast$ and
$R_a^*:\mathfrak{m}^\ast\to\mathfrak{m}^\ast$ denote the dual maps
of $L_a$ and $R_a$ respectively. If $v\in\mathfrak{m}_\ast$, then
$L_a^*(v)$ and $R_a^*(v)$ are $\sigma$-continuous functionals and
therefore, by the characterization of the predual
$\mathfrak{m}_\ast$  as the subspace of $\sigma$-continuous
functionals in $\mathfrak{m}^\ast$, it follows that $L_a^*(v),
R_a^*(v)\in\mathfrak{m}_\ast$.
Since $\operatorname{ad}_a=[a,\cdot]=L_a-R_a$ it follows that
$\operatorname{ad}_a^*=L_a^*-R_a^*$ and hence for any
$v \in \mathfrak{m}_\ast $, we have that $\operatorname{ad}_a^*(v)
= L_a^*(v)-R_a^*(v) \in \mathfrak{m}_\ast$.

This shows that  $\mathfrak{m}_\ast$ is a Banach Lie-Poisson
space with the Poisson bracket $\{f,g\}$ of $f,g\in
C^\infty(\mathfrak{m}_\ast)$ given by \eqref{general LP}. The
Hamiltonian vector field $X_f$ defined by the smooth function $f\in
C^\infty(\mathfrak{m}_\ast)$ is given by \eqref{general Hamiltonian vector field}.
\medskip

An {\bfi exact sequence of $W ^\ast$-algebras\/} is an exact
sequence of algebras in which all maps are $W ^\ast$-homomorphisms.

Let us analyze exact sequences of Banach Lie-Poisson spaces that
are preduals of $W ^\ast$-algebras. So assume that \eqref{lie
poisson sequence} is an exact sequence of Banach Lie-Poisson spaces
such that their  duals $\mathfrak{n} : = \mathfrak{c}^\ast$,
$\mathfrak{g}: = \mathfrak{b}^\ast$, $\mathfrak{h}: =
\mathfrak{a}^\ast$ are $W ^\ast$-algebras. By Theorem \ref{thm:
exact sequence of lie poisson}, the sequence
\eqref{dual lie algebra sequence},
where  $\iota: = p ^\ast$ and $\pi: =  j ^\ast$, is an exact
sequence of Banach Lie algebras in the category $\mathfrak{L}_0 $
but all objects in the sequence are
$W ^\ast$-algebras. This means that the maps $\iota$ and $\pi$ are
Banach Lie algebra homomorphisms and that $\iota
^\ast(\mathfrak{b}) \subset
\mathfrak{c}$, $\pi ^\ast(\mathfrak{a}) \subset \mathfrak{b}$. What
is not guaranteed, and is not true in general, is that the linear
continuous maps
$\iota$ and $\pi$ are homomorphisms of the
associative product structures of the $W ^\ast$-algebras
$\mathfrak{n}, \mathfrak{g}, \mathfrak{h}$. An example of a Lie
algebra homomorphism between $W ^\ast$-algebras that is not a
homomorphism for the associative product is given by $\phi:
\operatorname{gl}(n) \rightarrow \operatorname{gl}(n) $, $\phi(a) :
= \operatorname{tr}(a) \mathbb{I}$, where $\operatorname{gl}(n) $
is the algebra of $n \times n $ matrices and $\mathbb{I}$ the
identity matrix.

We shall assume now that $\iota$ and $\pi$ are
$\ast$-homomorphisms for the associative product structure. Then,
since they preserve the preduals, Proposition \ref{w star
homomorphism condition} insures that they are $W
^\ast$-homomorphisms. Conversely, assume that
\eqref{dual lie algebra sequence} is an exact sequence of $W
^\ast$-algebras. Then the maps $j $ and $p $ are homomorphisms of
Banach Lie algebras and, by Proposition \ref{w star homomorphism
condition}, their duals preserve the predual spaces, that is, this
is an exact sequence in the category $\mathfrak{L}_{0u}$.

Thus we are lead to consider exact sequences \eqref{dual lie algebra
sequence} of $W ^\ast$-algebras. Then $\ker \pi =\operatorname{im}
\iota $ is a $\sigma$-closed ideal in $\mathfrak{g}$ and thus there
exist a central projector $z \in \mathfrak{g}$ such that
$\operatorname{im} \iota = z \mathfrak{g}$ (Proposition 1.10.5 in
Sakai [1971]). The projector $1 - z $ is also central so that
denoting $(\operatorname{im} \iota) ^\perp : = (1-z) \mathfrak{g}$,
one has the direct sum splitting
\[
\mathfrak{g} = \operatorname{im} \iota \oplus (\operatorname{im}
\iota)^\perp
\]
into two $\sigma$-closed ideals of $\mathfrak{g}$.
It is easy to see that for any $x \in \operatorname{im} \iota$ and
$y \in (\operatorname{im}
\iota)^\perp $ we have $xy = 0 $.  Thus the $W ^\ast$-algebra
$\mathfrak{g}$ is the direct sum of two commuting $\sigma$-closed
ideals. In addition, the map $\pi|_{(\operatorname{im}
\iota)^\perp} : (\operatorname{im}
\iota)^\perp
\rightarrow \mathfrak{h}$ is an isomorphism of $W ^\ast$-algebras.
This proves the following.

\begin{proposition}
\label{extensions of w star algebras}
Any extension
\begin{equation*}
0 \longrightarrow \mathfrak{n} \stackrel{\iota}\longrightarrow
\mathfrak{g} \stackrel{\pi}\longrightarrow \mathfrak{h}
\longrightarrow 0
\end{equation*}
of the $W ^\ast$-algebra $\mathfrak{h}$ by the $W ^\ast$-algebra
$\mathfrak{n}$ (where $\iota$ and $\pi$ are $\ast$-homomorphisms of
the associative product structure) is isomorphic to the extension
\begin{equation*}
0 \longrightarrow \mathfrak{n} \stackrel{\iota}\longrightarrow
\mathfrak{n} \oplus \mathfrak{h} \stackrel{\pi}\longrightarrow
\mathfrak{h} \longrightarrow 0,
\end{equation*}
where $\iota(n) : = (n, 0)$ and $\pi(n, h): = h $, for $n \in
\mathfrak{n}$ and $h \in \mathfrak{h}$.
\end{proposition}

This proposition together with Theorems \ref{thm: exact sequence of
lie poisson} and  \ref{product theorem} immediately yields the
following result.

\begin{proposition}
\label{extensions of Lie-Poisson spaces}
Any extension
\[
0 \longrightarrow \mathfrak{a} \stackrel{j}\longrightarrow
\mathfrak{b} \stackrel{p}\longrightarrow \mathfrak{c}
\longrightarrow 0
\]
of the Banach Lie-Poisson space $\mathfrak{c}$ by the Banach
Lie-Poisson space $\mathfrak{a}$ such that the dual sequence is an
exact sequence of $W ^\ast$-algebras, is isomorphic to the extension
\[
0 \longrightarrow \mathfrak{a} \stackrel{j}\longrightarrow
\mathfrak{a} \oplus \mathfrak{c} \stackrel{p}\longrightarrow
\mathfrak{c}
\longrightarrow 0,
\]
where $j(a) := (a, 0)$ and $\pi(a,c) = c $, for $a\in
\mathfrak{a}$ and $c \in \mathfrak{c}$. This means that
$\mathfrak{b}$ is Poisson isomorphic to the product of the Banach
Lie-Poisson spaces $\mathfrak{a}$ and $\mathfrak{c}$ in the sense
of Theorem \ref{product theorem}.
\end{proposition}

\bigskip

\section{Extensions of Banach Lie algebras}
\label{Section: Extensions of Banach Lie algebras}

As we showed in the previous section, the
problem of extension of Banach Lie-Poisson spaces reduces to the
problem of extensions  of Banach Lie algebras in the subcategory
$\mathfrak{L}_0 $. We begin with some general
considerations in the category of Banach Lie algebras.

Let $\operatorname{aut}(\mathfrak{n}) : = \{D: \mathfrak{n} \rightarrow
\mathfrak{n}\mid D\text{~derivation~of~} \mathfrak{n}\}$ be the Banach Lie
algebra of all continuous linear derivations of $\mathfrak{n}$.
Recall that $D:  \mathfrak{n} \rightarrow \mathfrak{n}$ is a {\bfi
derivation\/} if
$D[\eta, \zeta] = [D\eta, \zeta] + [ \eta, D\zeta]$ for all $\eta, \zeta
\in \mathfrak{n}$. Denote by $\operatorname{int}(\mathfrak{n}): =
\{\operatorname{ad}_\eta \mid \eta \in \mathfrak{n}\}$ the
subalgebra of $\operatorname{aut}(\mathfrak{n}) $
consisting of {\bfi inner derivations\/}. In general, this is
not a closed subspace of $\operatorname{aut}(\mathfrak{n})$.
\textit{In this section we shall assume that
$\operatorname{int}(\mathfrak{n})$ is closed in
$\operatorname{aut}(\mathfrak{n})$} and hence
$\operatorname{int}(\mathfrak{n})$ is then a Banach Lie ideal
of $\operatorname{aut}(\mathfrak{n})$. Denote by
$\operatorname{out}(\mathfrak{n}) : =
\operatorname{aut}(\mathfrak{n})/\operatorname{int}(\mathfrak{n})$, the
Banach Lie algebra of {\bfi outer derivations\/} of $\mathfrak{n}$.
The norm on $\operatorname{int}(\mathfrak{n})$ and
$\operatorname{aut}(\mathfrak{n})$ is the usual operator norm induced
from the space $\operatorname{gl}(\mathfrak{n})$ of all linear
continuous maps of $\mathfrak{n}$ into itself.  The norm on
$\operatorname{out}(\mathfrak{n})$ is the quotient norm.

Let $C_{\mathfrak{g}}(\mathfrak{n}): = \{\xi \in \mathfrak{g}\mid [\xi,
\zeta] = 0 \text{~for~all~} \zeta\in \mathfrak{n}\}$ be the
{\bfi centralizer\/} of
$\mathfrak{n}$ in $\mathfrak{g}$ and $C(\mathfrak{n}):= \{\eta \in
\mathfrak{n}\mid [\eta, \zeta] = 0 \text{~for~all~} \zeta\in
\mathfrak{n}\} = C_\mathfrak{g}(\mathfrak{n}) \cap \mathfrak{n}$ be
the {\bfi center\/} of $\mathfrak{n}$; $C(\mathfrak{n})$ is a closed
ideal in  $C_{\mathfrak{g}}(\mathfrak{n})$ which is itself a Banach
Lie subalgebra of $\mathfrak{g}$. Consider the following
commutative diagram of exact sequences:

\unitlength=5mm
\begin{center}
\begin{picture}(19,10)
%first line
\put(1,9){\makebox(0,0){$0$}}
\put(9,9){\makebox(0,0){$0$}}
\put(16,9){\makebox(0,0){$0$}}
%second line
\put(-3,5){\makebox(0,0){$0$}}
\put(1,5){\makebox(0,0){$C ( \mathfrak{n})$}}
\put(8.7,5){\makebox(0,0){$C_{\mathfrak{g}}(\mathfrak{n})$}}
\put(16,5){\makebox(0,0){$C_{\mathfrak{g}}(\mathfrak{n})/C(\mathfrak{n})$}}
\put(22,5){\makebox(0,0){$0$}}
%third line
\put(-3,0){\makebox(0,0){$0$}}
\put(1,0){\makebox(0,0){$\mathfrak{n}$}}
\put(9,0){\makebox(0,0){$\mathfrak{g}$}}
\put(16,0){\makebox(0,0){$\mathfrak{h}$}}
\put(22,0){\makebox(0,0){$0$}}
%fourth line
\put(-3,-5){\makebox(0,0){$0$}}
\put(1,-5){\makebox(0,0){$\operatorname{int}(\mathfrak{n})$}}
\put(9,-5){\makebox(0,0){$\operatorname{aut}(\mathfrak{n})$}}
\put(16,-5){\makebox(0,0){$\operatorname{out}(\mathfrak{n})$}}
\put(22,-5){\makebox(0,0){$0$}}
%fifth line
\put(1,-9){\makebox(0,0){$0$}}
%first vertical arrows
\put(1,8){\vector(0,-1){2}}
\put(1,4){\vector(0,-1){3}}
\put(1,-1){\vector(0,-1){3}}
\put(1,-6){\vector(0,-1){2}}
%second vertical arrows
\put(9,8){\vector(0,-1){2}}
\put(9,4){\vector(0,-1){3}}
\put(9,-1){\vector(0,-1){3}}
%third vertical arrows
\put(16,8){\vector(0,-1){2}}
\put(16,4){\vector(0,-1){3}}
\put(16,-1){\vector(0,-1){3}}
%first horizontal arrows
\put(-2.4,5){\vector(1,0){2.3}}
\put(2.2,5){\vector(1,0){5.4}}
\put(10,5){\vector(1,0){3.7}}
\put(18.3,5){\vector(1,0){2.6}}
%second horizontal arrows
\put(-2.3,0){\vector(1,0){2.6}}
\put(2,0){\vector(1,0){6.4}}
\put(9.6,0){\vector(1,0){5.4}}
\put(16.6,0){\vector(1,0){4.4}}
%third horizontal arrows
\put(-2.4,-5){\vector(1,0){2.2}}
\put(2.2,-5){\vector(1,0){5.3}}
\put(10.2,-5){\vector(1,0){4.4}}
\put(17.3,-5){\vector(1,0){3.9}}
%horizontal labels
\put(5,.7){\makebox(0,0){$\iota$}}
\put(12,.7){\makebox(0,0){$\pi$}}
%vertical labels
\put(1.8,-2.4){\makebox(0,0){$\operatorname{ad}$}}
\put(10.1,-2.4){\makebox(0,0){$\operatorname{ad}|_\mathfrak{n}$}}
\put(16.8,-2.4){\makebox(0,0){$\sigma$}}
\end{picture}
\end{center}
\vspace{2in}
All maps that are not labeled are natural: they are inclusions
or projections on quotients. In particular, $\mathfrak{h}$ is
isomorphic as a Banach Lie algebra with the quotient
$\mathfrak{g}/\mathfrak{n}$ endowed with the quotient norm. The map
$\operatorname{ad}: \mathfrak{n}\rightarrow
\operatorname{int}(\mathfrak{n})$ is given by
$\eta \mapsto \operatorname{ad}_\eta : = [ \eta, \cdot ]$, for
$\eta \in \mathfrak{n}$. The map $\operatorname{ad}|_\mathfrak{n}:
\mathfrak{g}\rightarrow \operatorname{aut}(\mathfrak{g})$ is given
by $\xi \mapsto \operatorname{ad}_\xi|_\mathfrak{n}$ for $\xi
\in \mathfrak{g}$.

Finally, the map $\sigma : \mathfrak{h} \rightarrow
\operatorname{out}(\mathfrak{h})$ is defined in the following way. Let
$\eta \in \mathfrak{h}$ and $\xi \in \pi^{-1}(\eta) \subset \mathfrak{g}$.
Define $\sigma(\eta): = [\operatorname{ad}_\xi|_\mathfrak{n}] \in
\operatorname{out}(\mathfrak{h})$, where
$[\operatorname{ad}_\xi|_\mathfrak{n}]$ denotes the equivalence
class of the operator $\operatorname{ad}_\xi|_\mathfrak{n} \in
\operatorname{aut}(\mathfrak{g})$ with respect to the ideal
$\operatorname{int}(\mathfrak{n})$. The map $\sigma$ is well
defined for if $\xi_1, \xi_2 \in \pi ^{-1}(\eta)$, then exactness
of \eqref{dual lie algebra sequence} implies that $\xi_1 - \xi_2 \in
\mathfrak{n}$ and thus $ \operatorname{ad}_{\xi_1 -
\xi_2}|_\mathfrak{n} \in \operatorname{int}(\mathfrak{n})$.
Therefore, $[\operatorname{ad}_{\xi_1}|_\mathfrak{n}] =
[\operatorname{ad}_{\xi_2}|_\mathfrak{n}] $. The map $\sigma$ is
clearly a Lie algebra homomorphism. The defining equality for
$\sigma$, that is,  $(\sigma \circ \pi)(\xi) =
[\operatorname{ad}_\xi|_\mathfrak{n}]$ for every $\xi \in
\mathfrak{g}$, proves that the lower right square of the diagram is
commutative. Continuity of $\sigma$ is proved in the following way.
Since
$\mathfrak{g}$ is a Banach Lie algebra, there is a constant $C>0$
such that $\|[\xi,\zeta]\|\leq C\| \xi \| \| \zeta \|$, for all
$\xi, \zeta \in \mathfrak{g}$. Therefore, $\|\operatorname{ad}_\xi
\| \leq C \| \xi \| $, for all $\xi\in \mathfrak{g}$. Thus, for
$\eta\in \mathfrak{h}$ and $\xi \in \pi ^{-1}( \eta) $ arbitrary,
we have $\| \sigma( \eta) \| =
\|[\operatorname{ad}_\xi|_\mathfrak{n}] \| \leq \|
\operatorname{ad}_\xi|_\mathfrak{n}\| \leq \|
\operatorname{ad}_\xi \| \leq C \| \xi \|$, which shows, using the
isomorphism $\mathfrak{g}/\mathfrak{n} \cong \mathfrak{h}$ and the
definition  of the norm on the quotient, that $\sigma $ is
continuous.

Note that
$\operatorname{ad}|_\mathfrak{n}:
\mathfrak{g}\rightarrow \operatorname{aut}(\mathfrak{g})$ and
$\sigma: \mathfrak{h} \rightarrow \operatorname{out}(\mathfrak{h})$
are not surjective, in general.

Consider now a linear continuous section
$s: \mathfrak{h}\rightarrow \mathfrak{g}$, $\pi\circ s =
\operatorname{id}_{\mathfrak{h}}$, and assume that $\mathfrak{g} =
\mathfrak{n}\oplus s(\mathfrak{h})$, that is,
$s(\mathfrak{h})$ is closed and has as split complement the space
$\mathfrak{n}$. Define the isomorphism of Banach spaces $\psi:
\mathfrak{g}\rightarrow
\mathfrak{n}\oplus \mathfrak{h}$ by $\psi(\xi) = (\xi-
s(\pi(\xi)), \pi(\xi) )$,  whose inverse is given by
$\psi^{-1}(\zeta, \eta) :=
\zeta + s(\eta)$, for $\eta \in \mathfrak{h} $, $\zeta \in
\mathfrak{n}$, and $\xi\in \mathfrak{g}$.  Note that
$\psi^{-1}(0, \eta) = s(\eta)$ for any
$\eta \in \mathfrak{h}$ and $\psi^{-1}(\zeta, 0) = \zeta$ for any
$\zeta \in \mathfrak{n}$. Conversely, given an isomorphism of Banach
spaces $\psi: \mathfrak{g}\rightarrow \mathfrak{n}\oplus
\mathfrak{h}$ that is the identity on $\mathfrak{n}$, it determines
a linear continuous section $s: \mathfrak{h} \rightarrow
\mathfrak{g}$ by
$s(\eta): = \psi^{-1}(0, \eta)$ whose image $s(\mathfrak{h}) =
\psi^{-1}(\mathfrak{h})$ is a closed split subspace of
$\mathfrak{g}$ admitting $\mathfrak{n}$ as a complement. Thus
\textit{there is a bijective correspondence between the Banach
space isomorphisms $\psi: \mathfrak{g}\rightarrow
\mathfrak{n}\oplus \mathfrak{h}$ that are the identity on
$\mathfrak{n}$ and the linear continuous sections $s: \mathfrak{h}
\rightarrow \mathfrak{g}$ with closed split range admitting
$\mathfrak{n}$ as a complement.}

From this point on we shall assume that $\mathfrak{g}$ is
isomorphic to
$\mathfrak{n}\oplus\mathfrak{h}$ as a Banach space. Let us stress
that this sum is \textit{not\/} taken, in general, as a direct  sum
of Banach Lie algebras. The
isomorphism $\psi$ induces a Lie bracket on $\mathfrak{n}\oplus
\mathfrak{h}$ by
\begin{align}
\label{bracket on sum}
[(\zeta, \eta), ( \zeta', \eta')] : &= \psi\left( [ \psi^{-1}(\zeta,
\eta), \psi^{-1}(\zeta', \eta') ] \right)  \nonumber \\
&= \left( [ \zeta, \zeta']
 + \varphi(\eta)(\zeta') -\varphi(\eta')(\zeta) +
\omega(\eta, \eta'), [\eta, \eta'] \right),
\end{align}
where  $\omega : \mathfrak{h} \times \mathfrak{h} \rightarrow
\mathfrak{n}$ and $\varphi: \mathfrak{h} \rightarrow
\operatorname{aut}( \mathfrak{n})$ are defined by
\begin{equation}
\label{cocycle}
\omega( \eta, \eta') : = [s(\eta), s(\eta') ] - s([\eta, \eta'])
\end{equation}
\begin{equation}
\label{representation}
\varphi(\eta) : = [ s(\eta), \cdot]
\end{equation}
for any $\zeta, \zeta' \in \mathfrak{n}$ and $\eta, \eta' \in
\mathfrak{h}$.

Let us pose the inverse question. Given are two Banach Lie algebras
$\mathfrak{n}$ and $\mathfrak{h}$. Endow the direct sum Banach
space $\mathfrak{n}\oplus \mathfrak{h}$ with the continuous bilinear
skew symmetric operation given by \eqref{bracket on sum},
where $\omega: \mathfrak{h}\times \mathfrak{h} \rightarrow
\mathfrak{n}$ is a given continuous bilinear skew symmetric map and
$\varphi : \mathfrak{h}\rightarrow \operatorname{aut}(\mathfrak{n})$ is a
given continuous linear map. What are the conditions on $\omega$
and $\varphi$ so that this operation $[(\zeta, \eta), (\zeta',
\eta')]$ defined by the right hand side of \eqref{bracket on sum} is
a Lie bracket on
$\mathfrak{n}\oplus
\mathfrak{h}$? The answer to this question is given by the
following proposition.

\begin{proposition}
\label{proposition: Lie algebra extension}
Let $\mathfrak{n}$ and $\mathfrak{h}$ be two Banach Lie algebras,
$\omega: \mathfrak{h}\times \mathfrak{h} \rightarrow
\mathfrak{n}$ a continuous bilinear skew symmetric map, and
$\varphi : \mathfrak{h}\rightarrow
\operatorname{aut}(\mathfrak{n})$  a continuous linear map. Then
\begin{equation}
\label{bracket on sum abstract}
[(\zeta, \eta), ( \zeta', \eta')] = \left( [ \zeta, \zeta']
 + \varphi(\eta)(\zeta') -\varphi(\eta')(\zeta) +
\omega(\eta, \eta'), [\eta, \eta'] \right)
\end{equation}
for $\zeta \in \mathfrak{n}$ and $\eta \in \mathfrak{h}$ endows the
Banach space direct sum
$\mathfrak{g}:=\mathfrak{n}\oplus \mathfrak{h}$ with a Banach Lie
algebra structure if and only if
\begin{align}
\label{cocycle condition}
&\omega([\eta, \eta'], \eta'') + \omega([\eta', \eta''], \eta)
+ \omega([\eta'', \eta'], \eta') \nonumber \\
& \quad - \varphi(\eta)(\omega(\eta', \eta''))
 - \varphi(\eta')(\omega(\eta'', \eta))
 - \varphi(\eta'')(\omega(\eta, \eta'))  = 0
\end{align}
and
\begin{equation}
\label{representation condition}
\operatorname{ad}_{\omega(\eta, \eta')} + \varphi([\eta, \eta']) -
[\varphi(\eta), \varphi(\eta')] = 0
\end{equation}
for any $\eta, \eta', \eta'' \in \mathfrak{h}$. Consequently, the
Banach Lie algebra $\mathfrak{n}\oplus\mathfrak{h}$ is an extension
of the Banach Lie algebra $\mathfrak{n}$ by the Banach Lie algebra
$\mathfrak{h}$.
\end{proposition}

\noindent \textbf{Proof.} Since $\varphi$ is linear continuous and
$\omega$ is bilinear continuous, the bilinear skew symmetric
operation defined in \eqref{bracket  on sum abstract} is also
continuous.  So it is enough to show that \eqref{cocycle condition}
and \eqref{representation condition} are equivalent to the Jacobi
identity.

From the expression of the second component
in \eqref{bracket on sum abstract}, it follows that the Jacobi
identity gives no conditions on it since $\mathfrak{h}$ is a
Lie algebra. Thus only the first components need to calculated.  A
direct computation shows that the first component of $[[(\zeta,
\eta), (\zeta', \eta')], (\zeta'', \eta'')]$ equals
\begin{align}
\label{jacobi}
&[[\zeta, \zeta'], \zeta''] +
[\varphi(\eta)(\zeta'), \zeta'']
- [\varphi(\eta')(\zeta), \zeta'']
 + [\omega(\eta, \eta'), \zeta'']
\nonumber \\
&\qquad
 + \varphi([\eta, \eta'])(\zeta'')
- \varphi(\eta'')([ \zeta, \zeta'])
 - \varphi(\eta'')\varphi(\eta)( \zeta') \nonumber \\
&\qquad \qquad + \varphi(\eta'')\varphi(\eta')( \zeta)
  - \varphi(\eta'')(\omega(\eta, \eta'))
 +  \omega([\eta, \eta'], \eta'').
\end{align}
For $\zeta = \zeta' = \zeta'' = 0 $ this expression becomes
$\omega([\eta, \eta'], \eta'') - \varphi(\eta'')(\omega(\eta, \eta'))$.
Taking the circular permutations of \eqref{jacobi} and then
setting $\zeta = \zeta' = \zeta'' = 0 $ yields
\begin{align*}
&\omega([\eta, \eta'], \eta'') + \omega([\eta', \eta''], \eta)
+ \omega([\eta'', \eta'], \eta') \nonumber \\
& \quad - \varphi(\eta)(\omega(\eta', \eta''))
 - \varphi(\eta')(\omega(\eta'', \eta))
 - \varphi(\eta'')(\omega(\eta, \eta'))  = 0
\end{align*}
for any $\eta, \eta', \eta'' \in \mathfrak{h}$. This proves
\eqref{cocycle condition}.

In the sum of \eqref{jacobi} with the  two terms obtained from it
by circular permutations, there are expressions that add up to
zero. The sum of the first term in \eqref{jacobi} with its circular
permutations is zero since it is the Jacobi identity in the Lie
algebra $\mathfrak{n}$. By \eqref{cocycle condition}, the sum of the
ninth and the tenth term in \eqref{jacobi} plus their circular
permutations add also up to zero. Thus, the sum of
\eqref{jacobi} with its circular permutations equals
\begin{align*}
&[\omega(\eta, \eta'), \zeta''] + [\omega(\eta', \eta''), \zeta] +
[\omega(\eta'', \eta), \zeta'] \\
&  + \varphi([\eta, \eta'])(\zeta'')  + \varphi([\eta', \eta''])(\zeta)  +
\varphi([\eta'', \eta])(\zeta') \\
&- [\varphi(\eta), \varphi(\eta')](\zeta'') - [\varphi(\eta'), \varphi(\eta'')](\zeta) -
[\varphi(\eta''), \varphi(\eta)](\zeta') \\
& + [\varphi(\eta)(\zeta'), \zeta''] - [\varphi(\eta)(\zeta''), \zeta']
- \varphi(\eta)([ \zeta', \zeta''])
\\
&+[\varphi(\eta')(\zeta''), \zeta] - [\varphi(\eta')(\zeta), \zeta'']
- \varphi(\eta')([ \zeta'', \zeta])  \\
& + [\varphi(\eta'')(\zeta), \zeta'] - [\varphi(\eta'')(\zeta'), \zeta] -
\varphi(\eta'')([ \zeta, \zeta']).
\end{align*}
Each of the last three lines vanishes because $\varphi(\eta)$, $\varphi(\eta')$,
and $\varphi(\eta'')$ are derivations on $\mathfrak{n}$. Since $\zeta$,
$\zeta'$, and $\zeta''$ are arbitrary, from the remaining three top
lines we conclude that
\begin{equation*}
\operatorname{ad}_{\omega(\eta, \eta')} + \varphi([\eta, \eta']) -
[\varphi(\eta), \varphi(\eta')] = 0
\end{equation*}
for any $\eta, \eta', \eta'' \in \mathfrak{h}$, which proves
\eqref{representation condition}.

Conversely, suppose that $\omega: \mathfrak{h}\times \mathfrak{h}
\rightarrow \mathfrak{n}$ is a continuous bilinear skew symmetric
map and $\varphi : \mathfrak{h}\rightarrow
\operatorname{aut}(\mathfrak{n})$  is a continuous linear map
satisfying \eqref{cocycle condition} and \eqref{representation
condition}. A direct verification using
\eqref{cocycle condition} and \eqref{representation condition}
shows that the Jacobi identity  holds. Thus \eqref{bracket on
sum abstract} endows $\mathfrak{g}\oplus \mathfrak{h}$ with a
Banach Lie algebra structure.
\quad $\blacksquare$

\medskip

If $\mathfrak{n}$ is Abelian, then $\operatorname{ad}_{\omega(\eta,
\eta')} = 0$ and $\operatorname{aut}(\mathfrak{n}) =
\operatorname{gl}(\mathfrak{n})$, the
Banach space of all linear continuous maps from $\mathfrak{n}$ to
$\mathfrak{n}$. Thus the second condition becomes  $\varphi([\eta,
\eta']) = [\varphi(\eta), \varphi(\eta')]$ for all $\eta \in \mathfrak{n}$,
that is, $\varphi: \mathfrak{h} \rightarrow
\operatorname{gl}(\mathfrak{n})$ is a representation. The
first condition asserts that $\omega$ is a
$\mathfrak{h}$--cocycle relative to the representation $\varphi$.

We want to mention that such extenstions of Lie algebras were
considered in the purely algebraic context by Mori [1953], Shukla
[1966], and, more recently, by Alekseevsky, Michor, and Ruppert
[2000].

\section{Extensions of Banach Lie-Poisson spaces}
\label{Section: Extensions of Banach Lie-Poisson spaces}

Let us take Banach Lie-Poisson spaces $\mathfrak{a}$ and
$\mathfrak{c}$ and construct  the extension $\mathfrak{b}$ of
$\mathfrak{c}$ by $\mathfrak{a}$ in the sense of the category
$\mathfrak{B}$ of Banach Lie-Poisson spaces. We restrict our
considerations to the case $\mathfrak{b}$ equals the direct sum
$\mathfrak{c}\oplus\mathfrak{a}$ of Banach spaces.
Thus one has the Banach space exact  sequence
\begin{equation}
\label{exact sequence of blp}
0 \longrightarrow \mathfrak{a} \stackrel{j}\longrightarrow
\mathfrak{c} \oplus \mathfrak{a} \stackrel{p}\longrightarrow
\mathfrak{c} \longrightarrow 0,
\end{equation}
where $j(a) : = (0, a) $ and $p(c,a) : = c $.  The dual of this
sequence is
\begin{equation}
\label{exact sequence of bla}
0 \longrightarrow \mathfrak{n} \stackrel{p^\ast}\longrightarrow
\mathfrak{n} \oplus \mathfrak{h} \stackrel{j^\ast}\longrightarrow
\mathfrak{h} \longrightarrow 0,
\end{equation}
where $\mathfrak{h}:= \mathfrak{a}^\ast$, $ \mathfrak{n} =
\mathfrak{c}^\ast$, $p ^\ast(\zeta): = (\zeta, 0)$ and
$j^\ast(\zeta,\eta) = \eta$. Since $\mathfrak{a}$ and
$\mathfrak{c}$ are Banach Lie-Poisson spaces we have
\begin{equation}
\label{predual condition on blp}
\operatorname{ad}^\ast_{\mathfrak{h}} \mathfrak{a} \subset
\mathfrak{a} \qquad \text{and} \qquad
\operatorname{ad}^\ast_{\mathfrak{n}} \mathfrak{c} \subset
\mathfrak{c}.
\end{equation}
By Theorem \ref{thm: exact
sequence of lie poisson}, the question whether \eqref{exact
sequence of blp} is an exact sequence of Banach Lie-Poisson spaces
is equivalent to the question whether \eqref{exact sequence of bla}
is an exact sequence in the subcategory $\mathfrak{L}_0 \subset
\mathfrak{L}$. Proposition
\ref{proposition: Lie algebra extension} gives a necessary and
sufficient condition for \eqref{exact sequence of bla} to be an
exact sequence in the category $\mathfrak{L}$. The same sequence is
exact in the subcategory $\mathfrak{L}_0 $ if and only if $
\operatorname{ad}^\ast_{\mathfrak{n}\oplus\mathfrak{h}}
(\mathfrak{c}\oplus\mathfrak{a}) \subset
\mathfrak{c}\oplus\mathfrak{a}$. In order to see what this means we
use formula \eqref{bracket on sum abstract} to compute the
coadjoint representation on $\mathfrak{n}^\ast \oplus
\mathfrak{h}^\ast$ and get
\begin{equation}
\label{coadjoint rep}
\operatorname{ad} ^\ast_{(\zeta, \eta)}(c,a) = \left(
\operatorname{ad} ^\ast_\zeta c + \varphi(\eta) ^\ast c ,
\omega(\eta, \cdot ) ^\ast c - (\varphi(\cdot ) \zeta) ^\ast c +
\operatorname{ad} ^\ast_\eta a \right)
\end{equation}
for $c \in \mathfrak{n}^\ast$, $a \in \mathfrak{h}^\ast$, $\zeta\in
\mathfrak{n}$, and $\eta \in \mathfrak{h}$. The requirement that
this action preserves the preduals, together with the properties
\eqref{predual condition on blp} implies that
$\varphi(\eta)^\ast(\mathfrak{c}) \subset \mathfrak{c} $ and that
$\left(\omega(\eta, \cdot ) ^\ast - (\varphi(\cdot ) \zeta)
^\ast\right) (\mathfrak{c}) \subset \mathfrak{a}$ for all $\eta \in
\mathfrak{h}$ and all $\zeta \in \mathfrak{n}$. Taking here
alternatively $\eta = 0 $ and $\zeta = 0 $, the second condition
becomes $\omega(\eta, \cdot ) ^\ast(\mathfrak{c}) \subset
\mathfrak{a}$ and $(\varphi(\cdot ) \zeta)^\ast(\mathfrak{c})
\subset \mathfrak{a}$. We have proved the following theorem.

\begin{theorem}
\label{theorem: extensions of blp}
Given are two Lie-Poisson Banach spaces $\mathfrak{a}$ and
$\mathfrak{c}$ whose duals are the Banach Lie algebras
$\mathfrak{h}:= \mathfrak{a}^\ast$ and $\mathfrak{n}:=
\mathfrak{c}^\ast$ respectively, a continuous bilinear skew
symmetric map $\omega:
\mathfrak{h}\times \mathfrak{h} \rightarrow
\mathfrak{n}$, and a continuous linear map $\varphi :
\mathfrak{h} \rightarrow
\operatorname{aut}(\mathfrak{n})$ satisfying \eqref{cocycle
condition} and \eqref{representation condition}. The Banach space
$\mathfrak{c}\oplus\mathfrak{a}$ is an extension of $\mathfrak{c}$
by $\mathfrak{a}$ if and only if
\begin{equation}
\label{condition one}
\varphi(\eta)^\ast(\mathfrak{c})
\subset \mathfrak{c} , \qquad   (\varphi(\cdot )
\zeta)^\ast(\mathfrak{c})
\subset \mathfrak{a}, \qquad  \omega(\eta, \cdot )
^\ast(\mathfrak{c})
\subset \mathfrak{a}
\end{equation}
for all $\eta\in \mathfrak{h}=
\mathfrak{a}^\ast$ and
$\zeta \in \mathfrak{n} = \mathfrak{c}^\ast$.
\end{theorem}

Given $f \in C^\infty(\mathfrak{c}\oplus \mathfrak{a})$, define the
partial functional derivatives $\delta f/ \delta c \in
\mathfrak{c}^\ast$ and $\delta f/\delta a \in \mathfrak{a}^\ast$
by
\[
D_{\mathfrak{c}}f(c, a) (c') = \left\langle c', \frac{\delta
f}{\delta c}
\right\rangle \qquad \text{and} \qquad
D_{\mathfrak{a}}f(c, a) (a') = \left\langle a', \frac{\delta
f}{\delta a}
\right\rangle
\]
for all $c' \in \mathfrak{c}$ and all $a' \in \mathfrak{a}$, where
$D_{\mathfrak{c}}f(c, a)$ and $D_{\mathfrak{a}}f(c, a)$ denote the
partial Fr\'echet derivatives of $f $ at $(c, a) \in
\mathfrak{c}\oplus \mathfrak{a}$ respectively.

Theorem \ref{general theorem}, \eqref{bracket on sum abstract}, and
\eqref{coadjoint rep} thus immediately yield the following  theorem.

\begin{theorem}
\label{theorem: poisson bracket formula for restricted gl}
With the notations and hypotheses of Theorem \ref{theorem:
extensions of blp}, the Lie-Poisson bracket of $f, g \in C ^\infty
(\mathfrak{c}\oplus\mathfrak{a})$ is given by
\begin{align}
\label{lp bracket on extension}
&\{f, g\}(c, a) = \left\langle a, \left[ \frac{\delta f}{\delta a},
\frac{\delta g}{\delta a} \right] \right\rangle  \nonumber \\
& \qquad + \left\langle c, \left[ \frac{\delta f}{\delta c},
\frac{\delta g}{\delta c} \right] - \varphi\left(\frac{\delta
g}{\delta a}\right)\frac{\delta f}{\delta c}
+ \varphi\left(\frac{\delta
f}{\delta a}\right)\frac{\delta g}{\delta c}
+ \omega\left(\frac{\delta
f}{\delta a}, \frac{\delta g}{\delta a} \right)
\right\rangle
\end{align}
for $c \in \mathfrak{c}$ and $a \in \mathfrak{a}$. The Hamiltonian
vector field of $h \in C^\infty(\mathfrak{c}\oplus \mathfrak{a}) $
is given by
\begin{equation}
\label{lp Hamiltonian vector field}
X_h(c, a) =  -  \left(
\operatorname{ad} ^\ast_{\frac{\delta h}{\delta c}} c +
\varphi\left(\frac{\delta h}{\delta a}\right) ^\ast c , \,
\omega\left(\frac{\delta h}{\delta a}, \cdot \right) ^\ast c -
\left(\varphi(\cdot ) \frac{\delta h}{\delta c}\right) ^\ast c +
\operatorname{ad} ^\ast_{\frac{\delta h}{\delta a}} a
\right)
\end{equation}

\end{theorem}

\noindent \textbf{Example 1: Semidirect products of Banach
Lie-Poisson spaces.} Let us apply Theorem \ref{theorem: extensions
of blp} for the case $\omega= 0 $. Condition \eqref{cocycle
condition} is in this case vacuous and condition
\eqref{representation condition} asserts that $\varphi:
\mathfrak{h}\rightarrow \operatorname{aut}(\mathfrak{n})$ is a Lie
algebra homomorphism. One can define the {\bfi semidirect product
of $\mathfrak{h}$ with $\mathfrak{n}$\/} as the Banach Lie algebra
with underlying Banach space $\mathfrak{n} \oplus \mathfrak{h}$ and
bracket \eqref{bracket on sum abstract} with $\omega = 0$.

Denote, as before, by
$\mathfrak{a}$ and $\mathfrak{c}$ the predual spaces
of $\mathfrak{h}$ and $\mathfrak{n}$ respectively,
that is, $\mathfrak{h}:= \mathfrak{a}^\ast$ and $\mathfrak{n}:=
\mathfrak{c}^\ast$. In this case, only two conditions in Theorem
\ref{theorem: extensions of blp}
survive, namely,
\begin{equation}
\label{semidirect product conditions}
\varphi(\eta)^\ast(\mathfrak{c})
\subset \mathfrak{c}  \qquad \text{and} \qquad (\varphi(\cdot )
\zeta)^\ast(\mathfrak{c}) \subset \mathfrak{a}
\end{equation}
for all $\eta\in \mathfrak{h}= \mathfrak{a}^\ast$ and
$\zeta \in \mathfrak{n} = \mathfrak{c}^\ast$. The Banach
 space $\mathfrak{c} \oplus \mathfrak{a}$ is predual
to $\mathfrak{n}\oplus\mathfrak{h}$. It is a Banach Lie-Poisson
space relative to the Poisson bracket
\eqref{lp bracket on extension} and Hamiltonian vector field
formula \eqref{lp Hamiltonian vector field} with
$\omega = 0
$. This is the {\bfi semidirect product Banach Lie-Poisson space of
$\mathfrak{a}$ with $\mathfrak{c}$}.

An important particular case of this situation occurs when
$\mathfrak{n}$ is an Abelian Lie algebra, that is, $\varphi:
\mathfrak{h} \rightarrow L^\infty(\mathfrak{n})$ is a
Lie algebra representation; $L^\infty(\mathfrak{n})$
denotes the Banach algebra of all bounded linear operators on
$\mathfrak{n}$. In this case one obtains the semidirect product of
$\mathfrak{h}$ with the Banach space $\mathfrak{n}$ whose bracket
is given by
\[
[(\zeta, \eta), ( \zeta', \eta')] = \left(
\varphi(\eta)(\zeta') -\varphi(\eta')(\zeta), [\eta, \eta']
\right)
\]
for any $\zeta, \zeta' \in \mathfrak{n}$ and $\eta, \eta'
\in \mathfrak{h}$. If conditions \eqref{semidirect product
conditions} hold, one obtains the semidirect product Banach
Lie-Poisson space $\mathfrak{c} \oplus \mathfrak{a}$. The formula
for the Poisson bracket is \eqref{lp bracket on extension} with
$\omega= 0 $ and the first summand of the second term set also
equal to zero. The formula for the Hamiltonian vector field is
\eqref{lp Hamiltonian vector field} with $\omega= 0 $ and the first
term in the first component set equal to zero. These formulas
coincide with the ones found, for example, in Holm, Marsden, and
Ratiu [1998].

As a further special case, let us assume that $\mathfrak{h}$ is a
$W ^\ast$-algebra $\mathfrak{m}$ and that $\mathfrak{n}$ is a
Hilbert space $\mathcal{H}$. Additionally, let us fix a $W
^\ast$-representation $\varphi: \mathfrak{m} \rightarrow
L^\infty(\mathcal{H})$ of $\mathfrak{m}$ on the Hilbert space
$\mathcal{H}$. In this case, $\mathfrak{c}= \mathcal{H}$ and
$\mathfrak{a} = \mathfrak{m}_\ast $, where $\mathcal{H} \cong
\mathcal{H}_\ast \cong \mathcal{H}^\ast $ is equipped with the
trivial Poisson structure, since we consider $\mathcal{H}$ as an
Abelian Banach Lie algebra. Conditions \eqref{condition one} of
Theorem \ref{theorem: extensions of blp} reduce in this case  to
the single requirement
\[
(\varphi(\cdot)v) ^\ast \mathcal{H} \subset \mathfrak{m}_\ast
\qquad \text{for all} \qquad v \in \mathcal{H}.
\]
This condition can be expressed as follows: for any $v, w \in
\mathcal{H}$, there exists an element $b \in \mathfrak{m}_\ast $
such that
\[
\langle \varphi(x) v \mid w \rangle = \langle x , b \rangle,
\]
for any $x \in \mathfrak{m}$, where $\langle \cdot \mid \cdot
\rangle $ denotes the inner product on $\mathcal{H}$. But this
condition is satisfied since the representation
$\varphi$ is
$\sigma$-continuous (by definition) and thus the linear functional
$x \mapsto \langle
\varphi(x) v \mid w \rangle $ is $\sigma$-continuous too. This shows
that it is represented by an element $ b \in \mathfrak{m}_\ast $
(Sakai [1971]). Therefore we have constructed the semidirect
product Banach Lie-Poisson space $\mathcal{H} \oplus
\mathfrak{m}_\ast$.

Let us further specialize this situation for the case $\mathfrak{m}
= L^\infty(\mathcal{H})$. The predual space $\mathfrak{m}_\ast $ is
in this case
the Banach space of trace class operators $L^1(\mathcal{H})$ and
the duality pairing  between
$L^\infty(\mathcal{H}) $ and $L^1(\mathcal{H})$ is given by
$\operatorname{tr}(\rho X)$  for $\rho\in L^1(\mathcal{H})$ and $X
\in L^\infty(\mathcal{H}) $. Formula
\eqref{lp bracket on extension}  for $f, g \in C^\infty
(\mathcal{H}\oplus L^1(\mathcal{H}))$ becomes
\[
\{f, g \}(v, \rho) = \operatorname{tr} \left(\rho \left[
\frac{\delta f}{\delta
\rho}, \frac{\delta g}{\delta \rho}\right]\right) + \left\langle v\,
\Big| \,\frac{\delta f}{\delta \rho} \frac{\delta g}{\delta
v} -
\frac{\delta g}{\delta \rho} \frac{\delta f}{\delta v}
\right\rangle,
\]
where $\rho \in
L^1(\mathcal{H}) $ and $v \in \mathcal{H}$. Hamilton's equation
$\dot{f} =
\{f, h\}$ for the Hamiltonian $h \in
C^\infty(L^1(\mathcal{H}))$ can be equivalently written as the
system of equations
\[
|\dot{v}\rangle = -\left(\frac{\delta h}{\delta \rho}\right)^\ast
|v \rangle \qquad
\dot{\rho} = \left[ \frac{\delta h}{\delta \rho}, \rho \right] +
\Big| \frac{\delta h}{\delta \rho} \Big \rangle \langle v |.
\]

\bigskip

\noindent\textbf{Example 2: An extension of the restricted Banach
Lie-Poisson space.} Let $\mathcal{H}$ be a complex separable Hilbert
space endowed with a {\bfi polarization\/} (Pressley and
Segal [1986], Wurzbacher [1990]), that is, a direct sum decomposition
$\mathcal{H}= \mathcal{H}_+ \oplus \mathcal{H}_-$ into two closed
orthogonal subspaces. Denote by  $P_{\pm}: {\mathcal H} \rightarrow
{\mathcal H}_\pm$ the orthogonal projectors on $\mathcal{H}_\pm $;
hence
$P_+ + P_- = id$ and $P_+ P_- = P_- P_+= 0
$.  Denote by
$L^\infty({\mathcal H})$ and
$L^\infty({\mathcal H}_{\pm})$ the Banach Lie algebra of bounded
linear operators on
${\mathcal H}$ and ${\mathcal H}_{\pm}$ respectively, relative to
the commutator bracket. Let $L^2({\mathcal H})$ be the Banach Lie
algebra of linear Hilbert--Schmidt operators on
${\mathcal H}$, also relative to the commutator bracket.
Similarly, let $L^2({\mathcal H}_+, {\mathcal H}_-)$ and
$L^2({\mathcal H}_-, {\mathcal H}_+)$ be the Banach spaces of
Hilbert-Schmidt operators from ${\mathcal H}_+ $ to ${\mathcal
H}_-$ and ${\mathcal H}_- $ to ${\mathcal H}_+ $
respectively. Following Pressley and Segal [1986] we call
\begin{align}
\label{definition of h}
&\mathfrak{h}: = L^\infty({\mathcal H}, {\mathcal H}_+) := \{ X \in
L^\infty({\mathcal H}) \mid P_{\pm} X P_{\pm} \in L
^\infty({\mathcal H}_{\pm}), \nonumber  \\
&\qquad \qquad \qquad  P_+X P_- \in L^2({\mathcal H}_-,
{\mathcal H}_+), P_-X P_+ \in L ^2({\mathcal H}_+, {\mathcal H}_-)
\}
\end{align}
the
{\bfi restricted Banach Lie algebra\/}.
In this definition we write, for example, $P_+XP_-$ for
$P_+XP_-|_{{\mathcal H}_-}$ and similarly for the other terms. The
vector space $\mathfrak{h}$ is a Banach space relative to the norm
\begin{equation}
\label{definition of norm on h}
\|X\|: = \|P_+X P_+\|_\infty + \|P_-XP_-\|_\infty + \|P_+XP_-\|_2 +
\|P_-XP_+\|_2,
\end{equation}
where $\|\cdot \|_\infty $ and $\|\cdot \|_2 $ denote the operator
norm and the Hilbert--Schmidt norm in the various spaces. It is
easy to show that relative to the commutator bracket $[X, X']: =
XX' - X'X $, the space $\mathfrak{h}$ is  a Banach Lie algebra. It
also convenient to think of elements of $\mathfrak{h}$ as block
operators of the form
\begin{equation}
\label{matrix of operator x}
\left(
\begin{array}{cc}
X_+& X_{+-}\\
X_{-+}&  X_-
\end{array}
\right),
\end{equation}
where $X_{\pm}:= P_{\pm} X P_{\pm}\in L^\infty({\mathcal H}_{\pm})$,
$X_{+-} := P_+XP_- \in L^2({\mathcal H}_{-}, {\mathcal H}_+) $, and
$X_{-+}:= P_-X P_+ \in L^2({\mathcal H}_{+}, {\mathcal H}_-) $.

The Banach space $\mathfrak{n}: = L^1({\mathcal H}) $ of trace class
operators on ${\mathcal H}_+ $ endowed with the trace norm
$\|\cdot \|_1 $ and the negative of the commutator bracket $[ \rho,
\rho']:= -\rho \rho' + \rho' \rho $ is a also a Banach Lie algebra.
Define $\mathfrak{g}:= \mathfrak{n} \oplus \mathfrak{h}$ and
\begin{align}
\label{definition of phi}
\varphi: X \in \mathfrak{h} &\mapsto [P_+XP_+, \cdot ] \in
\operatorname{aut}(\mathfrak{n}) \\
\label{definition of omega}
\omega: (X,X') \in \mathfrak{h} \times \mathfrak{h} &\mapsto
P_+XP_- X' P_+ - P_+X'P_-X P_+ \in \mathfrak{n}.
\end{align}
The map $\varphi$ is linear and continuous and the map $\omega$ is
bilinear and continuous. These maps also satisfy the identities
\eqref{cocycle condition} and \eqref{representation condition}.
Indeed, to verify \eqref{representation condition}, for
arbitrary $\sigma\in L ^1 ({\mathcal H}_+)$ and $X, X' \in L
^\infty ({\mathcal H}, {\mathcal H}_+)$, taking into account that
the bracket operation on $\mathfrak{n}$ is the negative of the
commutator bracket, we have
\begin{align*}
&\operatorname{ad}_{\omega(X, X')}\sigma + \varphi([X, X']) \sigma-
[\varphi(X), \varphi(X')]\sigma \\
&\qquad  = - [\omega(X, X'), \sigma] +
[P_+[X,X']P_+, \sigma] - [[P_+XP_+, \cdot ], [P_+X'P_+, \cdot
]]\sigma \\
& \qquad  = -[P_+XP_- X' P_+ - P_+X'P_-X P_+, \sigma] \\
& \qquad \qquad + [P_+X(P_+
+ P_-) X' P_+ - P_+X'(P_+ + P_-)XP_-, \sigma] \\
& \qquad \qquad - [P_+XP_+, [P_+X'P_+, \sigma]] +
[P_+X'P_+, [P_+XP_+, \sigma]] \\
&\qquad = [P_+XP_+ X' P_+ - P_+X'P_+X P_+, \sigma] + [P_+XP_+, [
\sigma, P_+X'P_+]] \\
& \qquad \qquad + [P_+X'P_+, [P_+XP_+, \sigma]]\\
&\qquad = [ \sigma, [P_+X'P_+, P_+XP_+]] + [P_+XP_+, [
\sigma, P_+X'P_+]]
+ [P_+X'P_+, [P_+XP_+, \sigma]] = 0
\end{align*}
by the Jacobi identity.

To verify \eqref{cocycle condition} we compute separately the first
pair of terms for $X, X', X'' \in \mathfrak{h}= L^\infty({\mathcal
H}, {\mathcal H}_+)$ to get
\begin{align*}
&\omega([X, X'], X'') - \varphi(X)(\omega(X', X''))
=  P_+[X,X']P_- X'' P_+ - P_+X''P_-[X,X'] P_+ \\
&\qquad  \qquad - [P_+XP_+, P_+X'P_- X'' P_+ - P_+X''P_-X' P_+] \\
& \qquad = P_+X(P_+ + P_-)X'P_- X'' P_+
- P_+X'(P_+ + P_-)X P_- X'' P_+ \\
&\qquad \quad - P_+X''P_-X(P_+ + P_-)X' P_+
 + P_+X''P_-X'(P_+ + P_-)X P_+ \\
&\qquad \quad  - P_+XP_+X'P_-X''P_+ +  P_+XP_+X''P_-X'P_+\\
&\qquad \quad + P_+X'P_- X'' P_+XP_+ - P_+X''P_-X' P_+ X P_+\\
&\qquad = P_+X P_-X'P_- X'' P_+ - P_+X'P_+ X P_- X'' P_+
- P_+X' P_- X P_- X'' P_+ \\
&\qquad \quad - P_+X''P_-X P_+ X' P_+ - P_+X''P_-X P_-X' P_+
+ P_+X''P_-X' P_- X P_+ \\
&\qquad \quad +  P_+XP_+X''P_-X'P_+ + P_+X'P_- X'' P_+XP_+
\end{align*}
Rearrange the terms in the following manner:
\begin{align*}
&\omega([X, X'], X'') - \varphi(X)(\omega(X', X'')) =\\
& P_+X P_-X'P_- X'' P_+ + P_+X''P_-X' P_- X P_+
- P_+X' P_-X P_- X'' P_+ - P_+X''P_-X P_-X' P_+\\
& \qquad +  P_+XP_+X''P_-X'P_+ - P_+X'P_+ X P_- X'' P_+\\
& \qquad + P_+X'P_- X'' P_+XP_+ - P_+X''P_-X P_+ X' P_+\\
&\quad =P_+\big(X P_-X'P_- X'' + X''P_-X' P_- X
- X' P_-X P_- X''  - X''P_-X P_-X'  \big) P_+ \\
&\qquad + P_+\big( XP_+X''P_-X' - X'P_+ X P_- X'' \big) P_+\\
&\qquad + P_+\big( X'P_- X'' P_+X - X''P_-X P_+ X'\big)P_+.
\end{align*}
Adding the other two terms obtained by circular permutations  gives
zero; the summands cancel separately in the three groups emphasized
above. Thus Proposition \ref{proposition: Lie algebra extension} can
be applied thereby showing that $\mathfrak{g}$ is an extension of
$\mathfrak{h} = L^\infty({\mathcal H}, {\mathcal H}_+) $ by
$\mathfrak{n}= L^1({\mathcal H}_+)$. The Lie bracket on
$\mathfrak{n}\oplus \mathfrak{h}$ is given by  \eqref{bracket on sum
abstract} which in this case becomes
\begin{align}
\label{lie bracket for the restricted lie algebra}
[(X, \rho), ( X', \rho')] &= \left( -[ \rho, \rho']
 + \varphi(X)(\rho') -\varphi(X')(\rho) +
\omega(X, X'), [X, X'] \right) \nonumber\\
&= \big( -[ \rho, \rho'] + [P_+XP_+, \rho'] - [P_+X'P_+, \rho]
\nonumber \\
& \qquad \qquad + P_+XP_- X' P_+ - P_+X'P_-X P_+,\; [X, X'] \big)
\end{align}

The predual of $\mathfrak{h}$ is the Banach Lie-Poisson space
\begin{align}
\label{definition of predual of block lie algebra}
\mathfrak{a}:=L ^{1}({\mathcal H}, {\mathcal H}_+) : = \{\sigma \in&
L^\infty({\mathcal H}) \mid  P_{\pm} \sigma P_{\pm} \in
L^1({\mathcal H}_{\pm}), \nonumber \\
& P_+ \sigma P_- \in L^2({\mathcal H}_-, {\mathcal H}_+),
P_- \sigma P_+ \in L^2({\mathcal H}_+, {\mathcal H}_-)\}
\end{align}
relative to the pairing
\begin{equation}
\label{definition of the pairing on blocks}
\langle \sigma , X \rangle : = \operatorname{trace}(\sigma_+X_+ +
\sigma_-X_- + \sigma_{+-}X_{-+} + \sigma_{-+}X_{+-})
\end{equation}
for $X \in \mathfrak{h}$. The predual of $\mathfrak{n}$ is
$\mathfrak{c}: = {\mathcal K}({\mathcal H}_+)$, the Banach
Lie-Poisson space of compact operators on ${\mathcal H}_+ $. We
shall verify now the hypotheses of Theorem \ref{theorem: extensions
of blp}, that is,
\begin{itemize}
\item[{\rm (i)}] $\varphi(X)^\ast(\mathfrak{c}) \subset \mathfrak{c}
$, where $\varphi(X): L^1({\mathcal H}_+) \rightarrow  L^1({\mathcal
H}_+)$, so $\varphi(X) ^\ast: L^\infty({\mathcal H}_+) \rightarrow
L^\infty({\mathcal H}_+)$ and one needs to show
$\varphi(X)^\ast({\mathcal K}({\mathcal H}_+)) \subset {\mathcal
K}({\mathcal H}_+) $,
\item[{\rm (ii)}] $(\varphi(\cdot ) \rho)^\ast(\mathfrak{c})
\subset \mathfrak{a}$, where $\varphi(\cdot ) \rho :
L^\infty({\mathcal H},{\mathcal H}_+) \rightarrow L^1({\mathcal
H}_+)$, so $(\varphi(\cdot ) \rho) ^\ast :
L^\infty({\mathcal H}_+) \rightarrow (L^\infty({\mathcal H},
{\mathcal H}_+)) ^\ast$ and one needs to show  $(\varphi(\cdot )
\rho)^\ast({\mathcal K}({\mathcal H}_+) ) \subset L^1({\mathcal H},
{\mathcal H}_+)$
\item[{\rm (iii)}] $\omega(X, \cdot ) ^\ast(\mathfrak{c})
\subset \mathfrak{a}$, where $\omega(X, \cdot ) :
L^\infty({\mathcal H}, {\mathcal H}_+) \rightarrow L^1 ({\mathcal
H}_+) $, so $\omega(X, \cdot ) ^\ast :
L^\infty({\mathcal H}_+) \rightarrow (L^\infty({\mathcal H},
{\mathcal H}_+)) ^\ast$ and one needs to show $\omega(X, \cdot )
^\ast( {\mathcal K}({\mathcal H}_+)) \subset L^1({\mathcal H},
{\mathcal H}_+)$
\end{itemize}
for all $X\in \mathfrak{h}=
L^\infty({\mathcal H}, {\mathcal H}_+)$ and $\rho \in \mathfrak{n}
= L^1({\mathcal H}_+)$.

To verify (i) use the trace pairing, let
$X \in \mathfrak{h} $, $Y \in L^\infty ({\mathcal H}_+)$,  and
$\rho \in  L^1({\mathcal H}_+)$ to get
\[
\langle \varphi(X) ^\ast Y , \rho \rangle = \langle Y,
\varphi(X)\rho \rangle = \operatorname{trace}(Y[P_+XP_+, \rho]) =
\operatorname{trace}([Y, P_+XP_+] \rho)
\]
which shows that $\varphi(X) ^\ast Y = [Y, P_+XP_+] \in
L^\infty({\mathcal H}_+)$, that is,
$\varphi(X) ^\ast = - [ P_+XP_+, \cdot] $. Therefore, if $\varkappa
\in {\mathcal K}({\mathcal H}_+)$, we have $\varphi(X) ^\ast
\varkappa = - [ P_+XP_+, \varkappa] \in {\mathcal K}({\mathcal
H}_+)$, since
${\mathcal K}({\mathcal H}_+) $ is an ideal in $L^\infty({\mathcal
H}_+)$. This shows that  $\varphi(X)^\ast(\mathfrak{c}) \subset
\mathfrak{c} $.

To show (ii), use the same notations as before to get
\begin{align*}
\langle (\varphi(\cdot ) \rho)^\ast Y, X \rangle &= \langle Y,
\varphi(X) \rho\rangle = \operatorname{trace}(Y[P_+XP_+, \rho]
\rangle = \operatorname{trace}([\rho , Y]P_+ XP_+).
\end{align*}
Thus, if $\varkappa \in {\mathcal K}({\mathcal H}_+)$ we have
\[
\langle (\varphi(\cdot ) \rho)^\ast \varkappa, X \rangle
=\operatorname{trace}([\rho , \varkappa]P_+ XP_+)
\]
so, taking for $X $ operators of the form $X_-$, $X_{+-}$, and
$X_{-+}$, the right hand side of this relation vanishes. Therefore
\[
(\varphi(\cdot ) \rho)^\ast \varkappa =
\left(
\begin{array}{cc}
[\rho, \varkappa]& 0\\
0&  0
\end{array}
\right) \in L^1({\mathcal H}, {\mathcal H}_+)
\]
since $[\rho, \varkappa ] \in L^1({\mathcal H}_+)$.

Finally, to verify (iii), use the same notations as above,
let $X' \in L^\infty( {\mathcal H}, {\mathcal H}_+)$ to
get
\begin{align*}
\langle \omega(X, \cdot ) ^\ast Y, X' \rangle &= \langle Y,
\omega(X, X') \rangle = \operatorname{trace}(Y(P_+XP_- X' P_+ -
P_+X'P_-X P_+))\\
& = \operatorname{trace}(YX_{+-}X'_{-+} - X_{-+} Y X'_{+-}).
\end{align*}
Thus if $\varkappa \in {\mathcal K}({\mathcal H}_+)$, we have
\[
\langle \omega(X, \cdot ) ^\ast \varkappa, X' \rangle =
\operatorname{trace}(\varkappa X_{+-}X'_{-+} - X_{-+}
\varkappa X'_{+-})
\]
which shows that
\[
\omega(X, \cdot ) ^\ast \varkappa =
\left(
\begin{array}{cc}
0& \varkappa X_{+-}\\
- X_{-+} \varkappa &  0
\end{array}
\right) \in L^1({\mathcal H}, {\mathcal H}_+)
\]
because $\varkappa X_{+-} \in L^2({\mathcal H} _-, {\mathcal H}
_+)$, $X_{-+} \varkappa \in L^2({\mathcal H} _+, {\mathcal H} _-) $
since $X_{+-} \in L^2({\mathcal H} _-, {\mathcal H} _+)$, $X_{-+}
\in L^2 ({\mathcal H}_+, {\mathcal H}_-)$ and the $L^2$ operators
are an ideal in the algebra of bounded operators.

In view of Theorems \ref{theorem: extensions of blp} and
\ref{theorem: poisson bracket formula for restricted gl} the direct
sum $\mathfrak{c}\oplus \mathfrak{a} = {\mathcal K}({\mathcal H}_+)
\oplus L^1({\mathcal H}, {\mathcal H}_+)$ is a Banach Lie-Poisson
space relative to the bracket
\begin{align}
\label{lp bracket on restricted gl}
&\{f, g\}(\varkappa, \sigma) = \operatorname{trace}\left( \sigma
\left[
\frac{\delta f}{\delta \sigma},
\frac{\delta g}{\delta \sigma} \right] \right)  \nonumber \\
& \qquad + \operatorname{trace} \left( \varkappa \left( \left[
\frac{\delta f}{\delta
\varkappa},
\frac{\delta g}{\delta \varkappa} \right] -
\left[P_+\frac{\delta g}{\delta \sigma}P_+, \frac{\delta
f}{\delta \varkappa} \right] + \left[P_+ \frac{\delta
f}{\delta \sigma}P_+, \frac{\delta g}{\delta \varkappa} \right]
\right. \right. \nonumber \\
&\qquad \qquad \qquad \left. \left. + P_+ \frac{\delta
f}{\delta \sigma} P_- \frac{\delta g}{\delta \sigma} P_+ -
P_+ \frac{\delta g}{\delta \sigma} P_-  \frac{\delta
f}{\delta \sigma} P_+
\right) \right)
\end{align}
for $f, g \in C^\infty({\mathcal K}({\mathcal H}_+) \oplus
L^1({\mathcal H}, {\mathcal H}_+)) $ and  $\varkappa \in {\mathcal
K}({\mathcal H}_+)$, $\sigma \in L^1({\mathcal H}, {\mathcal H}_+)$.
The Hamiltonian vector field of
$h \in C^\infty({\mathcal K}({\mathcal H}_+) \oplus
L^1({\mathcal H}, {\mathcal H}_+)) $ is given by
\begin{align}
\label{lp Hamiltonian vector field on restricted gl}
X_h(\varkappa, &\sigma) =  -  \left(\left[ \varkappa,
\frac{\delta h}{\delta \varkappa} \right] + \left[ \varkappa,
P_+ \frac{\delta h}{\delta \sigma} P_+ \right] ,
\right. \nonumber \\
&\left.
\left(
\begin{array}{cc}
0& \varkappa (\delta h/ \delta \sigma)_{+-}\\
- (\delta h/ \delta \sigma)_{-+} \varkappa &  0
\end{array}
\right)
-
\left(
\begin{array}{cc}
[\delta h/ \delta \varkappa, \varkappa]& 0\\
0&  0
\end{array}
\right)
 + \left[ \sigma,
\frac{\delta h}{\delta \sigma}\right] \right).
\end{align}

\bigskip
\bigskip

\addcontentsline{toc}{section}{Acknowledgments}
\noindent\textbf{Acknowledgments.} We are grateful to Peter Michor
and Claude Roger for discussions concerning the subject of this paper
and some references. We want to thank the Erwin Schr\"odinger
Institute for Mathematical Physics for its hospitality. The second
author was partially supported by the European Commission and the
Swiss Federal Government through funding for the Research Training
Network
\emph{Mechanics and Symmetry in Europe} (MASIE) as well as the Swiss
National Science Foundation.

\bigskip

\section*{References}

\begin{description}

\item Alekseevsky, D., Michor, P.W., and Ruppert, W. [2000]
Extensions of Lie algebras, preprint, arXiv: math.DG/0005042 v2

\item Bratteli, O. and Robinson, D.W. [1979] \textit{
Operator Algebras and Quantum Statistical Mechanics I\/},
Springer-Verlag,  New York.

\item  Bratteli, O. and Robinson, D.W. [1981] \textit{
Operator Algebras and Quantum Statistical Mechanics II\/},
Springer-Verlag, New York.

\item de Azcarr\'aga, J.A. and Izquierdo, J.M. [1995] \textit{Lie
Groups, Lie Algebras, Cohomology and Some Applications\/},
Cambridge Monographs on Mathematical Physics, Cambridge University
Press.

\item Chernoff, P. R. and Marsden, J. E. [1974]
\textit{Properties of Infinite Dimensional Hamiltonian Systems},
Lecture Notes in Mathematics, {\bf 425}. Springer Verlag.

\item Emch, G. [1972] \textit{Algebraic Methods in Statistical
Mechanics\/}, Wiley Interscience.

\item Dubrovin, B.A., Krichever, I.M., and Novikov, S.P. [2001]
Integrable systems. I.
In \textit{Dynamical Systems IV}, pages 177--332,
Encyclopaedia Math. Sci., \textbf{4},
Springer, Berlin.

\item Holm, D.D., Marsden, J.E., and Ratiu, T.S. [1998] The
Euler-Poincar\'e equations and semidirect products, \textit{Adv. in
Math.\/}, \textbf{137}, 1--81.

\item Lie, S. [1890] \textit{Theorie der Transformationsgruppen,
Zweiter Abschnitt.} Teubner.

\item Marsden, J.E. and  Ratiu, T.S. [1994]  {\it
Introduction to Mechanics and Symmetry\/}, Texts in Applied
Mathematics, {\bf 17}, Second Edition, second printing 2003,
Springer-Verlag.

\item Mori, M. [1953] On the three-dimensional cohomology group of
Lie algebras, \textit{J. Math. Soc. Japan}, \textbf{5}, 171--183.

\item Odzijewicz, A. and Ratiu, T.S. [2003] Banach Lie-Poisson
spaces and reduction, \textit{Comm. Math. Phys.\/}, to appear.

\item Pressley, A. and Segal, G. [1986] \textit{Loop Groups\/}.
Oxford Mathematical Monographs, Clarendon Press, Oxford.

\item Sakai, S. [1971] {\it $C^\ast$-Algebras and
$W^\ast$-Algebras.\/} Ergebnisse der Mathematik und ihrer
Grenzgebiete, {\bf  60}, 1998 reprint of the 1971 edition, New York, NY:
Springer-Verlag.

\item Shukla, U. [1966] A cohomology for Lie algebras, \textit{J.
Math. Soc. Japan}, \textbf{18}, 275--289.

\item Takesaki, M. [1979] {\it Theory of Operator Algebras I},
Springer-Verlag.

\item Vaisman, I. [1994] {\it Lectures on the Geometry of
Poisson Manifolds}, Progress in Mathematics, {\bf 118}. Birkh\"auser
Verlag,
Basel, 1994.

\item Weinstein, A. [1983] The local structure of Poisson manifolds,
{\it Journ. Diff. Geom\/} {\bf 18}, 523--557.

\item Weinstein, A. [1998] Poisson geometry,
{\it Differential Geom. Appl.} {\bf 9}, 213--238.

\item Wurzbacher, T. [2001] Fermionic second quantization and the
geometry of the restricted Grassmannian, in \textit{Infinite
Dimensional K\"ahler Manifolds\/}, A. Huckleberry and T.
Wurzbacher eds., DMV Seminar \textbf{31}, Birkh\"auser, Basel

\end{description}

\end{document}